\documentclass[12pt]{amsart}
\usepackage{amsmath,amssymb,latexsym,graphicx}
\newtheorem{theorem}{Theorem}
\newtheorem{lemma}[theorem]{Lemma}

\def\qed{$\Box$}

\def\pf{{\bf Proof }}

\begin{document}
\title [Complex dimensions of real manifolds]{Complex dimensions of  real manifolds, attached analytic discs and parametric argument principle}

%thanks{The work of the second author was supported in part by the NSF Grants D%S-9971674,
%MS-0002195, and DMS-0604778. The work of the third author was supported in par% by
%he NSF Grants DMS-0200788 and DMS-0456868.} }
\author{Mark L.~Agranovsky} 

\address{Department of Mathematics, Bar-Ilan University, Rmat-Gan, 52900, Israel }
\maketitle
\begin{center}
Bar-Ilan University
\end{center}
\date{}

%\documentclass[12pt]{article}
%,a4paper, amsppt, reqno]{amsart}
%\documentstyle[12pt,]{article}
%\textwidth 6in
%\hoffset-.6in
%\renewcommand{\baselinestretch}{1.5}
%\usepackage{amsmath,amssymb,latexsym}
%\usepackage[dvips]{graphics}
%\input{Home}
%\textwidth 165mm
%\textheight 230mm
%\topmargin -5mm
%\evensidemargin -2mm
%\oddsidemargin -2mm
%\baselineskip24pt
%\input amssymb.stylatex
%\pagestyle{headings}
\newcommand{\bt}{\begin{Theorem}}
\newcommand{\et}{\end{Theorem}}
\newcommand{\bi}{\begin{itemize}}
\newcommand{\ei}{\end{itemize}}
\newcommand{\bea}{\begin{eqnarray}}
\newcommand{\eea}{\end{eqnarray}}
\newtheorem{Definition}{Definition}[section]
\newtheorem{Theorem}[Definition]{Theorem}
\newtheorem{Lemma}[Definition]{Lemma}%
\newtheorem{Exercise}{\sc Exercise}[section]
\newtheorem{Proposition}[Definition]{Proposition}
\newtheorem{Corollary}[Definition]{Corollary}
\newtheorem{Problem}[Definition]{Problem}
\newtheorem{Example}[Definition]{Example}
\newtheorem{Remark}[Definition]{Remark}
\newtheorem{Remarks}[Definition]{Remarks}
\newtheorem{Question}[Definition]{QuesMotition}
\newcommand{\be}{\begin{equation}}
\newcommand{\ee}{\end{equation}}

\vskip.25in

%\begin{center}
%\section{\protect\large\bf Introduction}
%\end{center}
\setcounter{equation}{0}
%\maketitle

% author two information

%\subjclass[2000]{Primary: 35L05, 44A12 Secondary: 35B05, 35S30}
%\date{}inva

\begin{abstract}

In the articles \cite{A1}, \cite{A2} characterizations of holomorphic and $CR$ functions were obtained in terms of 
analytic extensions from closed curves into analytic discs bounded by these curves. 
In particular, two open questions, the strip-problem and Globevnik-Stout conjecture 
were answered in real-analytic category.
The results from \cite{A1},\cite{A2}, reformulated as statements for  the graphs of the functions, give raise to
more general problem of deriving lower bounds for the complex dimensions ($CR$-dimensions) of real varieties of  complex manifolds, 
from properties of families of holomorphic 1-chains (analytic discs) which can be attached to the manifold. 
% More precisely,
%let $\Lambda$  be a real submanifold  of a complex manifold. Suppose that $\Lambda$ admits a  family 
%of attached (by their boundaries) holomorphic 1-chains (analytic discs).
%The question is: to what extent this family induces the complex structure on $\Lambda$?
We present a  solution of this  problem. 
As corollaries, we obtain solutions of some open problems and generalization of some known results on holomorphic and $CR$ 
functions and manifolds.
A surprising link between estimating of complex dimensions and a generalization of  
of the argument principle from a single holomorphic mapping to varieties of those is the core of the proof of the main results.

%Theorem 1 which  essentially says that a $CR-$ mapping of a bordered $CR$-manifold, which is homologically nontrivial
%on a "`real part"' of the manifold,  must induce nontrivial homomorphism of highest homologies of the boundary, unless
%the mapping is degenerate in the sense of the rank. This can be regarded as a generalization, for parametric families of holomorphic mappings, 
%of the argument principle in the following its version: a holomorphic mapping $f$ of the unit disc $\Delta \subset \mathbb C$ to $\mathbb C,$ 
%smooth in $\overline \Delta,$ is constant, if $f\vert_{\partial \Delta}$ has  topoloigical degree zero. The very relation of the problem
%under study with topological form of argument principle seems to be  of indpendent interest. 

\end{abstract}
\maketitle

\vskip 1cm

\section{Introduction}\label{sec:1}

In this article we establish and study the link between the three following subjects:

\begin{itemize}

\item
Characterization of complex or $CR$ manifolds and holomorphic or $CR$ functions by vanishing moments on families of closed curves.

\item Global properties of smooth families of analytic 
discs attached to real submanifolds of complex manifolds.

\item
A version of the generalized argument principle for smooth families
of holomorphic mappings (parametric argument principle).
\end{itemize}

Albeit the third subject seems to be  of special and independent interest, in this 
publication we view it rather as the main tool, while  
the first two subjects will be in the main focus.

\medskip
\noindent{\it Motivation.}
This article develops the ideas of the articles \cite{A1} and \cite{A2} of the author. 
There, the following two question about holomorphic functions and their boundary values were answered in real analytic category.  We discuss these two question in more details in the next sections, hence  in this introductory section we will touch on them  briefly. 

The first question, called "strip-problem", (see  \cite{AV},
\cite{AG},\
\cite{E2},\cite{E1},
\cite{G1},\cite{JG}, \cite{G3}, \cite{G4},\cite{G6},\cite{G7}, 
\cite{T1},\cite{T2}, \cite{Z1}), asks whether a function of one complex variable is analytic if it extends as analytic function
from one-parameter family of Jordan curves? The condition of analytic extendibility from the curves  can be formulated 
as condition of vanishing of all complex moments hence the question is about a  version of Morera theorem.

The second question, Globevnik-Stout conjecture, or, more generally, boundary Morera problem,
(see \cite{AV}, \cite{ABC}, \cite{AS}, \cite{AS1},
\cite {GS}, \cite {GS1}, \cite{KM}), \cite{D1}, \cite{D2}, \cite{D3}, \cite{NR}, \cite{R}, \cite{St2}), 
is: if a function on the boundary of a domain in 
$\mathbb C^n$, or more generally, on a $CR$-manifold, is the boundary value of a holomorphic function (is $CR$-function, i.e. satisfies the boundary 
$\overline \partial$-equation) if it extends
holomorphically (satisfies zero moment condition) from a family of closed curves which are boundaries of analytic discs in $\mathbb C^n$ ?
(Such analytic discs are called attached.)
The Globevnik-Stout conjecture, which generalizes an earlier result for the complex ball, due to Nagel and Rudin, is a special case of the boundary Morera problem when the attached analytic discs  are obtained by complex linear sections.

In this article, we generalize the results of \cite{A1}, \cite{A2} and obtain new results in this direction by looking at  the problems in a   
more general context. The scheme is as follows.

\medskip
{\it Step 1. From functions to manifolds.}  
The first step is based on the simple idea that properties of functions often can be reworded as those of surfaces, if one passes from functions to their graphs. 
In particular,  the question whether a function is holomorphic  translates as whether its graph is a complex manifold, i.e. the tangent spaces are complex linear subspaces in $\mathbb C^n$? Correspondingly, asking  whether a function  is $CR-$ function is equivalent to asking
whether its graph has the dimension of complex tangent spaces ($CR$-dimension) not less than the $CR$-dimension of the manifold where the function is originally defined.

Taking into account  that the graphs of the analytic extensions into attached analytic discs, are themselves analytic  discs,
attached to the graph of our function,  we immediately arrive to a general problem of estimating $CR$-dimension of real manifolds in $\mathbb C^n$ in terms of families of attached analytic discs.  Information about $CR$-dimension is  important in many problems, for instance, in the problem of
characterization of boundaries of complex varieties (see  \cite{DH0}, \cite{DH}, \cite{D1}, \cite{D3}), \cite{Ha}, \cite{HL}, \cite{W}) which we also consider in this article. 

\medskip
\noindent
{\it Step 2. From manifolds to mappings.}
Observe that the
$CR$-dimension of a real manifold  at its point can be "`measured"' by the dimension of the $\mathbb C$-linear span of the tangent space at this point:
the more is $CR$-dimension, the less is the dimension of the span. Therefore, the $CR$-dimension is a "`measure" of linear dependence of the tangent space over $\mathbb C$, i.e. measure of the degeneracy of our real manifold over $\mathbb C$.  

Now, assume that our real manifold is covered by the boundaries of analytic discs constituting a smooth family parametrized by point in a real manifold of parameters. Each analytic disc in this family is a holomorphic image
of the standard complex unit disc and hence we deal with a smooth family of holomorphic mappings of the unit disc into $\mathbb C^n$. 
This smooth family of holomorphic mappings can be regarded as a single mapping holomorphically depending on one complex variable in the unit disc 
and smoothly-on the rest of  the real variables (the parameters of the family).

Our manifold is covered by the images of the unit circles, therefore it can be viewed as the image, under the above mapping, 
of the Cartesian product of two manifolds: the unit circle and the  manifold of parameters. 
Thus, according to the earlier remark, the $CR$-dimension of our real manifold is related to the $\mathbb C$-degeneracy of the tangent spaces, which in turn is related to the rank of the above partially holomorphic mapping.

\medskip

\noindent
{\it Step 3. Parametric argument principle.}
Thus, we have linked the question of estimating of $CR$-dimension of a real manifold, equipped by  attached analytic discs, 
with the question about degeneracy, in the sense of the rank, of a certain mapping.
The  third step  relates the  latter question with  a generalized  argument principle for holomorphic mappings.

This relation is the key ingredient of our considerations and requires more detailed explanation. 
In fact, what we mean is the following corollary (a weak version) of the argument principle: 
if  a holomorphic mapping of the closed unit disc into $\mathbb C^n$
maps the unit circle into a homologically trivial loop then the mapping is constant (has rank zero). 

Indeed, otherwise
the image of the unit disc is a 1-dimensional complex variety and hence 
has the positive intersection index with any other complex variety of  the complementary
dimension, intersecting the image of the unit disc but disjoint from the image its boundary. 
However, this contradicts to the general form of the argument principle. Indeed, it says that 
the above intersection index equals to the linking coefficient between
the complementary variety and the loop-the image of the unit circle. 
But this loop is homologically trivial and therefore the linking coefficient, and hence the interesection index, must  0,
which is the  contradiction. 

In our situation we deal with a family (manifold) of holomorphic mappings rather than with a single mapping. 
This manifold of mappings can be regarded also as
as a single mapping but defined on  the Cartesian product of the unit disc and the manifold of parameters. 
The mapping depends holomorphically 
on the first, complex, variable in the unit disc, and  smoothly - on the rest of the (real) variables-parameters.  

Our key result is a generalization of the above version of the argument principle, from a single mapping to  families of those,
by proving the following fact.
Suppose that the above mapping is homologically trivial
on the boundary  which is the Cartesian product of the circle and of the parameterizing manifold. Suppose that, however,
the trajectory of a fixed point in the unit disc, 
when the parameters run the parameterizing manifold, sweeps up a cycle which is not homological to 0
in the entire image. Then  the entire mapping degenerates in the sense of its rank. 
In particular, the dimension of the image falls down. 

In less formal terms, the result says the following: if
a smooth parametric family of holomorphic mappings of the unit disc into $\mathbb C^n$ collapses on the unit circle, meaning that
the images of the unit circle sweep up a variety with trivial highest homologies, 
but this homological degeneracy  comes not from 
the parameters, then the entire family of mappings degenerates in the sense of the rank. In this case the images of the unit disc sweep up
a variety of the dimension less than maximally possible.  
Notice that the ordinary argument principle corresponds to the parameterizing manifold consisting of a single point.

This result, in view of the above reduction, implies estimating of $CR$-dimensions and, correspondingly, leads to
solution of the above mentioned related problems.

\section{Preliminary discussion}

\subsection{Description of the main problem}

Let $\Lambda$ be a real submanifold of a complex Stein manifold $X=X^n$ of the complex dimension $n.$  
$CR$-structure on $\Lambda$ is defined by the maximal complex subspaces
of the tangent space of $\Lambda$. The manifold $\Lambda$ is $CR$ is the dimension of complex tangent subspaces is the same at all points
and this dimensiond is called $CR$-dimension. 
$CR$-functions on a $CR$-manifold are smooth functions satisfying the tangential $CR$-equations in the complex tangent spaces. 
For continuous functions, the $CR$-conditions can be formulated in a weak sense, in terms of orthogonality to $\overline \partial$-closed forms.

If $\Lambda$ is a complex submanifold then $CR$-functions are the same as holomorphic functions. 
By Bochner-Severi theorem $CR$-functions on the  boundary of a domain in $\mathbb C^n$
coincide with the boundary values of holomorphic functions in this domain.

Analytic disc $D=\varphi(D)$ in $X$ is the image of the closed unit complex disc $\overline \Delta$
under holomorphic immersion $\varphi:\overline \Delta \mapsto X.$ The disc is attached to a real manifold $\Lambda$
if $\partial D=\varphi(\partial \Delta) \subset \Lambda.$ 
Analytic discs are very powerful tool for study $CR$-manifolds, and the literature on this subject is very
large.

The main problem, we study in this article, can be formulated as follows:

\noindent
{\bf Problem.} {\it Given a real submanifold of a complex manifold, equipped by a family of attached analytic discs,
derive lower bound for its $CR$-dimension.}

Let us discuss this problem in more details.
First of all, when the manifold is the graph of a function, then the question turns to a characterization of $CR$-functions.
If $n=1$  then both the manifold and the attached analytic discs lie in the same
complex plane and we arrive to the problem of characterization of holomorphic functions in terms of analytic extensions from
closed curves.  

When the manifold is a boundary of a domain in $\mathbb C^n$ then we deal with the Morera type problem of
characterization of boundary values of holomorphic functions.

In the authors's articles \cite{A1},\cite{A2} such characterization were obtained under assumption of real analyticity
and, in particular, two open questions, the strip-problem and Globevnik-Stout conjecture,
were answered in real analytic category. These problems, which we have mentioned already in Introduction,  
concern characterization of holomorphic functions in planar domains and boundary values of holomorphic  functions
in domains in $\mathbb C^n$ (see \cite{A1},\cite{A2}, \cite{AV},
\cite{ABC}, \cite{AS1},\cite{AS}, \cite{MLA}, \cite{AG},\cite{BTZ},\cite{D1},
\cite{E2},\cite{E1},
\cite{G1}-\cite{G7}, 
\cite{GS},\cite{GS1}, \cite{KM}, \cite{T1},\cite{T2}, \cite{Z1}).

The case when the manifold $\Lambda$ is the graph of a function $f$ is of special interest.
If $f$ analytically extends to an attached analytic disc, then the graph of the
extension is an analytic disc attached to the graph of $f$.  
Further, the function $f$ is holomorphic or $CR$ means that its graph is a complex or $CR$-manifold.   
Therefore, if one passes from the language of functions to the language of the manifold, then one is led
to a problem of determining the $CR$-dimension  of a real manifold admitting  families of attached analytic 
discs with certain properties.

Let us make another important remark. The minimal number of parameters of a family $\mathcal F$ of curves covering a manifold $\Lambda$
is $\dim \mathcal F=\dim \Lambda-1.$
If one is allowed to shrink analytic discs to a point then their tangent complex vectors converge to a tangent complex vector
in the tangent space of the manifold at the limit point. Having at hands rich enough family of discs one can induce complex tangent spaces
of arbitrary possible dimension.
However, to induce a  tangent complex line, by the above limit process,  one needs at least $(\dim \Lambda+1)$-parametric family of attached discs.
On the contrary, we use the  families containing no infinitesimally small discs. This circumstance makes the problem much more difficult.

Our basic result is Theorem 2 describing families of attached analytic discs which enforces the $CR$-dimension to be positive.
This result makes sense for $d \leq n$ since otherwise the $CR$-dimension is automatically positive. 
For the range of the dimension $d$, we describe those families $\mathcal F$ of attached analytic discs
depending on minimally possible number of parameters:

$$  dim \mathcal F= \dim \Lambda-1 \ \mbox{ or} \   \dim \mathcal F=\dim \Lambda$$ 

\noindent
which imply existing of complex lines in the tangent spaces. 
Then we apply this result to  the case of general dimensions $d$ and of
arbitrary $CR$-dimensions.

\subsection{Examples.}

The existence of families, even large ones, of attached analytic discs does not guarantee nontrivial $CR$-structure, 
induced by the ambient complex manifold.

Here is a  simple example.
Define the 2-dimensional manifold in $\mathbb C^2$ as
$$\Lambda=\{(z_1,|z_1|) \in \mathbb C^2 : \ r \leq |z_1| \leq  R \},$$ 
where $r$ and $R$ are given positive numbers. 
The analytic disc $$D_t=\{(z,t) \in \mathbb C^2: |z| \leq t\},$$
is attached to $\Lambda$ because $\partial D_t =\{|z|=t\} \subset \Lambda.$ These analytic discs constitute the 1-parameter smooth family. 
Nevertheless, the manifold $\Lambda$  is totally real, i.e. has zero $C$-dimension.

In fact, the manifold $\Lambda$, which geometrically is a truncated  cone (see Fig.1, on the left), 
is the graph of the function $|z|$ over the annulus $\{r \leq |z| \leq R\}$ in $\mathbb C$ and the trivial observation is that
$|z|$ extends analytically from any circle $|z|=t$ but is not analytic.
The analytic discs attached to  $\Lambda$  
are obtained by by parallel sections $Re z_2=t$ of  the solid 3-dimensional cone $\{Im z_2=0, |z_1| \leq Re z_2\}.$

Moreover, even essentially larger families analytic discs may not detect complex structure. 
For instance, it was observed in \cite{G1} that the function $f(z)=1/\overline z$ extends analytically inside any circle surrounding 0.
These circles constitute 3-parameter family. Nevertheless the function $f$ is not analytic. 
The geometric interpretation  is as follows. The graph $$\Lambda=\{(z,1/\overline z) \in \mathbb C^2: r \leq |z| \leq R\}$$ of the function
$1/\overline z$
admits 3-parameter family of attached analytic discs - the graphs of the corresponding analytic extensions.  Nevertheless, $\Lambda$
is not a complex manifold, the opposite, it is totally real.
 
Let us give more examples to illustrate the main problem and the main result of the article. 

\medskip
\noindent
{\bf Example 1.}
Consider the 2-dimensional torus $\Lambda=T^2$ in $\mathbb C^2$ realized as the 
Shilov boundary of the bidisc $\Delta^2$, where $\Delta$ is 
the unit disc in the complex plane, i.e. $T^2=\partial \Delta \times \partial \Delta$.  The 2-dimensional  manifold $T^2$
admits 1-parameter family $S^1 \times \Delta$ of  attached analytic discs but is totally real, $\dim_{CR}T^2=0.$ 

\medskip
\noindent
{\bf Example 2.}
Let $\Lambda$ be an embedded 2-dimensional annulus
$$\Lambda=\{(z,z) \in \mathbb C^2 : 1 <|z|<3 \},$$
and let
$$D_t=\{(z,z) \in \mathbb C^2: |z-2t|=1 \},$$  
where $|t|=1.$

Thus, $\Lambda$ is just the graph of the function $z$ over the annulus, and $D_t$ are the graphs of discs inscribed in the annulus. 
As in Example 1, the 1-parameter family of analytic discs $D_t$ is parametrized by points $t \in S^1$, and each analytic disc
$D_t$ is attached (in fact, belongs) to $\Lambda$. The manifold $\Lambda$ is complex, i.e. $\dim_{CR}\Lambda.$ 
It is illustrated on Fig.1, on the right.

In the next two examples the manifolds are 3-dimensional.
\medskip

\noindent
{\bf Example 3.}
Consider the 3-manifold in $\mathbb C^3$:
$$\Lambda=\{(z_1,z_2,|z_1|^2): |z_1|^2+|z_2|^2=1\}.$$ 
This manifold admits 3-parameter family of attached analytic discs $D_t$, parametrized by the points 
$$t=(t_1,t_2) \in S^3=\{(t_1,t_2) \in \mathbb C^2:|t_1|^2 + |t_2|^2=1\},$$ namely:
$$D_t=\{(\lambda_1\zeta, \lambda_2 \zeta, |\lambda_1|^2) \in \mathbb C^3: |\zeta| <1\},$$
where 
$\lambda_j=t_j/\sqrt{|t_1|^2 + |t_2|^2}, \ \ j=1,2.$ 

The manifold $\Lambda$ is the graph of the function $|z_1|^2$ over the unit sphere $S^3$ in $\mathbb C^2$. This function is constant 
on the sections of $S^3$ by 1-dimensional complex linear subspaces and corresponding attached analytic discs $D_t$ are just the graphs of the constant extensions
into those linear sections. $\Lambda$ is totally real, $\dim_{CR}\Lambda=0,$ at each point except the circle $(0,z_2,0), \ |z_2|=1.$

\medskip
\noindent
{\bf Example 4.} Define 3-manifold in $\mathbb C^3$ by
$$\Lambda=\{(z_1,z_2,z_1) \in \mathbb C^3: |z_1|^2+|z_2|^2=1\}.$$
Take $t=(t_1,t_2) \in S^3 \subset \mathbb C^2$ and 
let $L_t$ be the (unique) complex line tangent to the 3-sphere 
$$\{|z_1|^2+|z_2|^2=(1/4)(|t_1|^2+|t_2|^2)\}$$
at the point $(1/2)(t_1,t_2).$
Define
$$D_t=\{(z_1,z_2,z_1): (z_1,z_2) \in L_t, \ |z_1|^2+|z_2|^2 \leq \sqrt{3}/2\}.$$
Then $\{D_t\}_{t \in S^3}$ is a 3-parametric family of analytic discs attached to $\Lambda.$ 
The manifold $\Lambda$ is the graph of the function $z_1$ over $S^3.$ The analytic discs $D_t$ are graphs
of $z_1$ over the sections of the unit ball in $\mathbb C^2$ by complex lines tangent to the sphere $\frac{1}{2} S^3.$  
The manifold $\Lambda$ is maximally complex, $\dim_{CR}\Lambda=1.$

\medskip
Let us briefly analyse the examples. First, notice that in all of them the manifold  $\Lambda$ and the families of the attached discs  are real-analytic.
In  Examples 1 and 3 the manifolds are totally real, while  in Examples 2 and 4 the $CR$-dimensions are positive. 

Now, let us  examine the dimensions.
In  Example 1 and  Example 2 the dimension of manifolds $\Lambda$ is $d=2$ and  the dimension of the family of discs is $k=d-1=1$. 
The boundaries of the attached  discs are 1-dimensional curves and  the families of these curves depend on one real parameter. 
The curves cover  2-dimensional manifolds and therefore no dimensional degeneracy happens in both examples. 
Nevertheless, $\Lambda$ in Example 1 is totally real while in Example 2 it is a complex manifold.

However, an important  observation is that, albeit in  both Examples 1 and 2 the families of the discs are parametrized by the same 
closed manifold $S^1,$ in Example 1 the manifold $\Lambda$ is closed 
while in Example 2 it is not. In particular, the {\it Brouwer degree}  of the mapping 
$$\Phi:S^1 \times S^1 \mapsto \Lambda,$$
which parametrizes the family of the boundaries $D_t,$ equals $\deg \Phi=1$ in Example 1 ($\Phi$ is identical mapping in this case)  and equals 
$\deg \Phi=0$ in Example 2, because 
the topological degree of mappings from closed compact manifold to compact manifold with nonempty boundary is 0.

Now turn to Examples 3 and 4. Here in both cases the manifold of parameters is 3-dimensional sphere $S^3.$
The dimension of the manifolds $\Lambda$ is in both cases  $d=3$ and  the dimension of the families is  
$k=3.$ Contrary to Examples 1 and 2, here  we face dimensional
degeneracy because  $\dim (S^1 \times S^3)=k+1=4 > \dim \Lambda=d=3.$  
However, the results again are different: in Example 3 $\Lambda$ is totally real (except the circle 
$\{(0,z_2,0):|z_2|=1\}$), while in Example 4  $\Lambda$ is maximally complex (the $CR$-dimension equals 1).
Hence, there must be other reasons , besides degeneracy,  which affect on $CR$-dimensions. 

To understand the crucial difference between Examples 3 and 4, let us  us look at them from the topological point of view.
Observe that in Example 3 the projections of the discs $D_t$ on the plane $z_1,z_2$ are all complex  lines through 0 and
their union coincides with the unit ball in $\mathbb C^2.$   Therefore, the discs $D_t$ sweep up in $\mathbb C^3$ a 
5-manifold $\tilde \Lambda \subset \mathbb C^3,$ which is contractible and therefore  topologically trivial. In particular,
the centers of the discs $D_t$ fill the segment $[0,1],$ representing a topologically trivial loop.

The opposite, in Example 4, the projections of the discs on $\mathbb C^2$ are complex lines tangent
to the sphere of radius $1/2.$  The 4-dimensional manifold $\tilde \Lambda \subset \mathbb C^3$ in this example is noncontractible, as it contains a hole corresponding
the ball $|t| <1/2.$ Moreover, the centers of the discs $D_t$
form a 3-cycle which bounds no 4-cycle in $\tilde \Lambda$, because this 3-cycle surrounds the 4-dimensional hole.

%The illustration on Fig.1 and Fig. 2 shows the manifolds $\Lambda$ which are graphs of functions $|z|$ and $z$ corespondinlgy, over an annulus 
%$\{r \leq |z| \leq R \}$ in the plane. 
%Two types of families of attached analytic discs, homologically trivial and homologically nontrivial ones, are shown on the pictures. 
%The first family is $D_t=\{(z,t):|z| \leq t\}$, the second one is described in Example 2.

\begin{figure}[h]
\centering
\scalebox{0.5}{\includegraphics{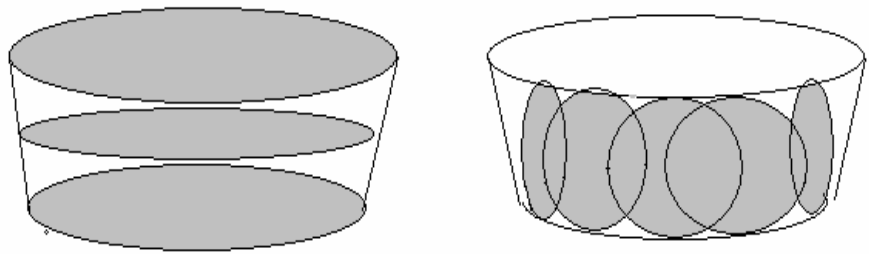}}
\label{fig:cones}
\caption{Homologically trivial (left) and homologically nontrivial (right) families of attached analytic discs}
\end{figure}

\subsection{Comments on the results}

Thus, we observe in the above four examples that the manifold $\Lambda$ carries no complex structure  in the cases  when
the family of the attached discs 
is either homologically nontrivial but non-degenerate on $\Lambda$ (Example 1), or degenerate on $\Lambda$ but homologically trivial (Example 3). 

Remind, that by degeneracy on $\Lambda$ we understand any of
the two following situations. The first one is the dimensional degeneracy, when $k+1 >d$, as in Examples 3 and 4, 
where $k$ is number of parameters of the family. 
The second type is when $k+1=d$ but the Brouwer degree is $0$, as in Example 2. Both types of degeneracy
are unified by saying that the $(k+1)$-th homology group of $\Lambda$ vanishes. 

On the other hand, in Example 2 and Example 4, where  the manifolds $\Lambda$ has positive $CR$-dimension, 
the families of attached analytic discs are both {\it degenerate and homologically nontrivial}.

It appears that just the last two properties are crucial and the main result of this article, Theorem 2, says 
that families of attached discs, possessing these two properties, do induce  $CR$-structure. 
The dimension of this structure depends on the size and the structure of the family of the attached analytic discs.

One should mention that the phenomenon is essentially global, as the property of homological nontriviality is a 
global property (though, the condition of degeneracy is obviously local). Moreover, examples  
show that one can not expect local characterizations of the considered type.

%\subsection{Parametric argument principle.}
The results on attached analytic discs rest on Theorem 1  which is the key theorem in  the article. 
Theorem 1 can be interpreted as a  {\bf parametric argument principle}. 

This relation was mentioned in Introduction, but let us explain it now in a more detailed way. 
The argument principle for complex manifolds can be understood as follows:
if $V_1^s,V_2^{n-s}$ are two transversally intersecting bordered complex submanifolds in a complex manifold $X^n$
then then the linking number of the boundaries equals to the number of the intersection points:  
$c(\partial V_1,\partial V_2)=\#(V_1 \cap V_2),$ counting multiplicities. This reflects the fact that the intersection index
of two complex varieties is nonnegative.

In particular,  if $\varphi:\Delta \mapsto X$ is a  nondegenerate holomorphic mapping
of  the unit complex disc to a complex manifold, smooth up to the boundary, 
then the linking number of the curve $\Gamma=\varphi(S^1)$
and the boundary of any disjoint from $\Gamma$  complex variety  is nonnegative. 
(Alexander and Wermer proved in \cite{AW} that this condition is
also necessary for $\Gamma$ to bound a 1-dimensional complex manifold or chain.) 

Therefore, if $\Gamma$ is homologically trivial, i.e. $H_1(\Gamma;\mathbb Z)=0,$ so that all the above linking numbers are $0,$ then
no complex chain is bounded by $\Gamma$ and hence $\varphi$ degenerates: $\varphi(\overline \Delta)=\Gamma.$ 
This fact can be regarded as a weak version of the argument principle for holomorphic mappings.

Now, suppose that,  instead of a single mapping, we have a smooth parametric family $\varphi_t(\zeta)$ of analytic mappings from $\Delta$ to
$\mathbb C^n$,  parametrized by points
 $t$ from a real manifold $M.$ Such objects naturally appear as parameterizations of smooth families of attached analytic discs.  
Theorem 1 generalizes
the above degeneracy phenomenon, from the case of single analytic mapping, $\dim M=0$, to multi-parametric families, $\dim M >0.$.

Loosely speaking, the collapse (homological triviality) of the union $\cup_{t \in M}\varphi_t(S^1)$ of the images of the unit circle 
may have two sources:  the dependence on the parameters $t$ and/or the holomorphic dependence on the complex variable $\zeta$. 
Therefore, to expect the generalization of the above weak form of the argument principle, 
one has to exclude the possible affect of the real parameters $t$.
The condition of {\it homological nontriviality in $t$} of the family just serves this purpose.  
Namely,  the condition requires that for some (and then for any) fixed $z_0 \in \Delta$ the mapping 
$$M \ni t \mapsto \varphi_t(z_0) \in X $$ 
maps $M$ to a nontrivial cycle in the union of the images,  $\cup_{t \in M} \varphi_t(\overline \Delta).$ 

It turns out that in this situation collapse of the boundary can be caused only by dimensional degeneracy of the mapping 
$(t,\zeta) \mapsto \varphi_t(\zeta).$ 
More precisely, our key result,Theorem 1,
claims that {\it if the family $\varphi_t, \ t \in M,$ of analytic mappings of the unit disc, homologically degenerates on the boundary, $S^1 \times M$, and is homologically nontrivial in $t \in M$, then it degenerates, in the sense of rank, in the solid manifold $\Delta \times M.$}

Theorem 1, applied to the parameterizing mappings  of the families of attached analytic discs, delivers the needed information about $CR$-dimension of the target manifold
because the degeneracy which claims Theorem 1 implies nonzero complex tangent spaces in the target manifold $\Lambda=\cup_{t \in M} \varphi_t(S^1).$

\medskip

\subsection{ Content of the article}

The article is organized as follows. 
Section 2 is devoted to theorems for manifolds. All objects, manifolds and  families of the analytic discs
are assumed real-analytic.  

All theorems in Section 2 are derived from Theorem 1 which is the basic fact behind all the results of this article.
Theorem 1, which we prove in Section 4,  says that if a $CR$-mapping
from a Cartesian product $\Delta \times M,$ where is the unit complex disc $\Delta$ and $M$ is a real compact manifold, 
degenerates on the boundary $S^1 \times M$
and is homologically nontrivial then the degeneracy propagates inside the Cartesian product. 

Theorem 1 generalizes, from a single analytic function to parametric families of analytic functions, 
the following corollary of the argument principle: if an analytic function in the unit disc, smooth up to its closure,
degenerates on the boundary (i.e., either maps it to  a point, or has degree 0) then it is constant.

Theorem 2 presents a geometric version of Theorem 1 and says that  if $d \leq n$ then real $d$-manifold in $n$-dimensional complex space  everywhere  the positive $CR$-dimension (i.e., is nowhere totally real)
if it admits a degenerate homologically nontrivial family of attached analytic discs. 
Corollary \ref{C:Morera_curves} is a partial of Theorem 2 for $n=2.$  It gives a characterization of complex curves in $\mathbb C^2$
as a real 2-manifolds admitting nontrivial families of attached analytic discs.

Theorem 2 stops working when $d >n.$ Indeed the $CR$-dimension is always at least $d-n$ and therefore it is automatically positive
when $d>n$. 
That is why we prove Theorem 3 which extends Theorem 2 for the range $d \geq n,$ in an effecitve way.  Namely, Theorem 3
derives, from the existence of certain families of attached analytic discs, the lower bound $\dim_{CR} \Lambda > q$ which  is, of course,
nontrivial  for $q > d-n.$

Theorem 4 is a special case of Theorem 3, for the case $\dim \Lambda=2p-1, \ q=p-1$. It describes
maximally complex manifolds, which, due to known results of Wermer \cite{W}, Harvey and Lawson \cite{HL} are
exactly boundaries of complex $p$-chains. 
Therefore, Theorem 4 characterizes those real manifolds which are boundaries of complex manifolds, in terms of attached analytic discs. 
Results of this type , for the case of analytic discs, obtained by complex linear sections, are earlier obtained 
by Dolbeault and Henkin \cite{DH0}, \cite{DH} and T.-C. Dinh \cite{D1},
\cite{D2},\cite{D3}. 

Section 3 is devoted to characterization of $C$R-functions on $CR$-manifolds. The obtained results are just special
cases of corresponding theorems from Section 2,  when the manifolds under consideration are graphs of functions.
We give characterizations $CR$-functions and boundary values of holomorphic functions, in terms of analytic extension in families of
attached analytic disc (one-dimensional extension property).
In particular, these results  contain solutions of the strip-problem and of Globevnik-Stout conjecture, for real-analytic functions (see
\cite{A1},\cite {A2}). 
Similar characterization of $CR$-functions
were obtained by Tumanov in \cite{T3}. Albeit, in \cite{T3} pretty special families of attached discs (thin discs) were exploited.

In Section 4 we prove Theorem 1 which serves  the base for all other results. We present two versions of the proof. The first one
works for the general case, while the second one-for closed families of discs. Both proofs 
exploit the duality between homologies and cohomologies, but the first proof uses Poincare duality and intersection indices and is more geometric,
while the second one uses de Rham duality and integrals of differential forms, computing linking numbers, and is more analytic. 

Both versions of the proof develop the ideas and constructions from
the earlier articles \cite{A1}, \cite{A2} of the author. However, here we not only generalize results from \cite {A1} and \cite{A2}
to higher dimensions, and from $CR$-functions to $CR$-manifolds, 
but we also  modify the approach. In particular, we tried to expose the constructions and the arguments, as much as possible, 
in a coordinate free form.  

As we have mentioned, Theorem 1 is our key tool for estimation $CR$-dimensions by means of attached analytic discs. 
In the concluding section we discuss the generalization of Theorem 1 to a version of argument principle
for $CR$-manifolds and $CR$-mappings.

\bigskip
\section{Basic notions and notations}\label{sec:2}

%\section{The key theorem. Parametric argument principle}
In the notations for manifolds, we will use upper index for the dimensions. For real manifolds it will mean the real dimension, and
and for complex manifolds - the complex one. However, sometimes we will be omitting this index. 

Let $X=X^n$ be a complex Stein manifold of the complex dimension $n.$ For simplicity we assume $X$ realized as a submanifold of the Euclidean space, $X \subset \mathbb C^N.$  
Let $\Lambda=\Lambda^d \subset X$ be a real smooth  $d-$dimensional submanifold, maybe with the boundary.

For each  point $b \in \Lambda$ we denote $T_b \Lambda$ the real tangent space to $\Lambda$ at the point $b$ and by 
$$T_b^{\mathbb C}\Lambda= T_b(\Lambda) \cap i T_b^{\mathbb C}(\Lambda)$$- 
the maximal complex tangent subspace.

Denote $$c(b)=c_{\Lambda}(b)= \dim_{\mathbb C} T_b^{\mathbb C}(\Lambda)$$
the complex dimension of the complex tangent space. We call $c(b)$ the $CR$-dimension of $\Lambda$ at the point $b$.It is easy to see that
$$c(b) \ge d-n.$$
For $d \geq n$, the dimension $c(b)$ takes the minimal possible value $c(b)=d-n$ if and only if $T_b\Lambda$ generates maximally possible complex linear space:
$$T_b\Lambda + i T_b\Lambda= T_b X.$$
In this case the manifold $\Lambda$ is called {\it generic.} When $d=n$ , the generic manifolds  
have no nontrivial complex tangent spaces. Such manifolds are called {\it totally real}.
For $ d >n$ the manifold $\Lambda^d$ is never totally real, while when $d \leq n$ the manifold may be totally real and may be  not.

When $c(b)=const$ then $\Lambda$ is called $CR$-manifold. In this case, a smooth function $f$ on $\Lambda$ is called $CR$-function
if it satisfies the tangential $CR$-equation,
$$\overline Zf=0$$
for any local tangent field $Z \in T^{\mathbb C}\Lambda.$

\medskip
We  will be dealing with families of smooth attached analytic discs, $D_t$, smoothly parametrized by points $t$ of a real-analytic
manifold, $M=M^k.$ The manifold  $M^k$ will be assumed oriented and compact, with boundary or without it.
Each analytic disc $D_t, \ t \in M$ is parametrized by a holomorphic immersion
$$\varphi_t:\Delta \mapsto D_t,$$
which is smooth up to the boundary of the unit disc.
Smooth  dependence  of the family $D_t$
on the parameter $t$ means that the mapping 
$$\Phi : \overline \Delta \times M \mapsto X,$$ 
defined by
$$\Phi(\zeta,t)=\varphi_t(\zeta)$$
belongs to $C^{r}(\overline \Delta \times M).$
The mapping $\Phi$ is a {\it parameterization} of the family $D_t$.

The order $r$ of smoothness can be $r=\infty$ or $r=\omega$. In the latter case we deal with real-analytic family.   
Throughout the paper the families of attached analytic disc will be real-analytic.

Let $\Lambda^d \subset X^n$ be a real submanifold as above. The fact that the discs 
$$\{D_t\}_{t \in M}$$ are attached to $\Lambda$ simply means that
$\Phi(\partial \Delta \times M) \subset \Lambda$. We assume more, namely, that $\Lambda$ is covered by the closed curves $\partial D_t$ i.e.
$$\Phi(\partial \Delta \times M)=\Lambda.$$ 
In this case we say that $\Lambda$ {\it admits} the family 
$$\mathcal F= \mathcal F_M= \{D_t\}_{t \in M}.$$

The mapping $\Phi$ is assumed {\it regular}  on $S^1 \times M^k$
meaning that
$$rank \ d\Phi\vert_{S^1 \times M^k}(u) = \dim \Lambda=d$$  
for $u \in \Phi^{-1}(\Lambda \setminus \partial \Lambda).$

From the construction, the dimension $k=\dim M$ of the family $\mathcal F$  cannot be less than $d-1$:
$$ k \ge d-1.$$  By the regularity condition, the dimensions of the fibers $\Phi^{-1}(b), b \in \Phi(\overline \Delta \times M^k)$  
do not exceed 1.

We will use the following notations:
\begin{equation}\label{E:notation}
\Sigma =\Delta \times M, \  b\Sigma=S^1  \times M,  b_0 \Sigma=\overline \Delta \times \partial M, \partial (b\Sigma)=S^1 \times \partial M.
\end{equation}

\section{The key theorem. Parametric argument principle}

The main results of this article follow from the following theorem about propagation of boundary degeneracy, which we regard
as parametric argument principle (see Introduction).

Everywhere in the sequel, the homology groups are understood with coefficients in $\mathbb Z$.

\begin{theorem} Let $M=M^k$ be a compact oriented real-analytic manifold and
$\Phi:\overline \Delta \times M^k \mapsto X^n, \ k  \le  n,$ be a real-analytic mapping, holomorphic in the variable 
$\zeta \in \overline \Delta.$ Assume that $H_{k-2}(M^k, \partial M^k)=0.$ Suppose that
$\Phi$ maps $b\Sigma=S^1 \times M^k$ onto a real-analytic variety $\Lambda^d, \ \ d \leq n,$   and $\Phi$ is regular on $b \Sigma$. 

%Denote $$\Sigma =\Delta \times M, \  b\Sigma=S^1  \times M^k, \ b_0 \Sigma=\overline \Delta \times \partial M^k,
%\Lambda=\Phi(S^1 \times M^k).$$ The manifold $b\Sigma$ itself may have the boundary $\partial (b\Sigma)=S^1 \tmes \partial M^k.$

Suppose also that 

\noindent
1. $\Phi$ homologically degenerates (is homologically trivial) on $b\Sigma$ in the following sense: the induced homomorphism of 
the relative homology groups:

$$\Phi_*:H_{k+1}(b\Sigma, \partial (b\Sigma) )
 \mapsto H_{k+1}(\Lambda,  \Phi(\partial(b\Sigma))$$
is zero, $\Phi_*=0.$

2. $\Phi$ is homologically nontrivial on $\Sigma,$ meaning that the 
induced homomorphism 
$$\Phi_* : H_k(\overline \Sigma, b_0 \Sigma) \mapsto 
H_k(\Phi(\overline \Sigma),\Phi(b_0\Sigma))$$
is not zero, $\Phi_* \neq 0.$
Here we have used the notations introduced in (\ref{E:notation}).

Then $rank_{\mathbb C} d \Phi\vert_{\overline \Sigma} < rank_{\mathbb R} \ d\Phi\vert_{b\Sigma}.$
%$\Phi$ degenerates in $\Delta \times M^k$ i.e. $rank \ d\Phi < k+2$ everywhere in $\overline \Delta \times M^k.$
None of the  conditions 1 and 2 can be omitted. 

\end{theorem}

\medskip
\begin{Remark}\label{R:condition_1}
The condition 1 is equivalent to the following: either 

\noindent
a) either $\Phi$ degenerates on $b\Sigma$ dimensionally, i.e. 
$rank \ \Phi\vert_{b\Sigma} < k+1$ and correspondingly $d=\dim \Lambda < k+1$ 

\noindent
or 

\noindent
b) $rank \ \Phi\vert_{b\Sigma} =k+1$ everywhere on $b\Sigma$ except the critical set, $\dim \Lambda=k+1$, but $\Phi$ degenerates on
$b\Sigma$ topologically meaning that the topological degree $deg \ \Phi\vert_{b\Sigma}=0.$
\end{Remark}

\pf
The $(k+1)$-th homology group of the $(k+1)$-manifold $b\Sigma=S^1 \times M^k$ is generated by the fundamental homology class
$$[b\Sigma]=[S^1 \times M^k].$$
The maximal rank of $\Phi$ on $b\Sigma$ is $k+1$ and $\Lambda$ is at most $(k+1)$-dimensional.
If $\Lambda$ has dimension $k$ or  less, then the $(k+1)$-th homologies of $\Lambda$ are trivial. Otherwise,
$\Lambda$ has dimension $k+1$ and then 
$$\Phi_*[b\Sigma]=\deg \Phi \cdot [\Lambda],$$ 
where $\deg \Phi$ is the Brouwer 
degree. Therefore, in this case  the condition 1 means $\deg \Phi=0.$  
\qed

%Taking into account Remark \ref{R:condition_1}, let us present an equivalent form of Theorem 1, in terms of families of holomorphic mappings fo the %union disc (parametric argument principle). 

%\bigskip
%{\bf Theorem 1'.}   
%{\it Let  
%$$\varphi_t \mapsto \overline \Delta \mapsto X^n, \ t \in M^k $$ be a real-analytic family
%of holomorphic mappings of the closed unit disc, parameterized by the manifold $M=M^k$ as in Theorem 1.
%Suppose that $\Lambda=\cup _{t \in M} \varphi(\partial \Delta)$ is a manifold of the dimension $d \leq n.$
%Suppose that either $d < k+1$ 
%or $d=k+1$ and the mapping $\partial D \times M \ni (\zeta,t) \mapsto \varphi_t(\zeta) \in \Lambda$ has the Brouwer degree 0.
%Suppose also that the family if homologically nontrivial, meaning that if $\zeta_0 \in \overline \Delta$ is a fixed point then
%the trajectory $\varphi_t(\zeta_0), \ t \in M,$ sweeps up the cycle which is not homological to 0 in the union 
%$\cup_{t \in M} \varphi(\overline \Delta)$ relatively to the same kind of union when $t$ runs $\partial M$. 
%Then 
%$$dim_{\mathbb C} span_{\mathbb C} (\partial_{\zeta}\varphi_t(\zeta), \partial_{t_1}\varphi_t(\zeta),
%\cdots, \partial_{t_k}\varphi_t(\zeta)) < d, \ (\zeta,t) \in \overline \Delta \times M,$$
%where $t_1,\cdots,t_k$ are local parameters on $M^k.$}

\begin{Remark}
If we take $k=0$ in Theorem 1  
then $M$ is just a point,  $\Sigma=\Delta$, and $\Phi$ is an analytic function $\varphi$ in the unit complex disc $\overline \Delta.$ Also, $\Lambda=\Phi(\partial \Delta)$ and  $d=dim \Lambda$ is either $d=0$ or $d=1.$ Condition 1
in this case means that either $\Lambda$ is a point, or $\Lambda$ is a curve with nonempty boundary so that the mapping $\Phi:S^1 \mapsto \Lambda$
has the topological degree 0. 
The conclusion of our theorems  says that $rank_{\mathbb C} d \Phi < d$ and therefore $rank_{\mathbb C} d\Phi=0$
which means that $\Phi=const,$ in accordance with the classical argument principle.

%The condition 1 says that either $\Phi(\partial \Delta)$ is a point, and then, of course, $\Phi=const$
%or $\Phi(\partial \Delta)$ is a curve with boundary, and $\Phi$ had Brouwer degree $0$ on $\partial \Delta.$
%Then Theorem 1 concludes that $rank \Phi <2$ in $\Delta$. 
%Then the image of $\Delta$ is 1-dimensional and then $\Phi=const,$ by the argument principle for analytic functions.

Thus, Theorem 1 can be viewed as a generalization of the argument principle from a single analytic function $\Phi(\zeta)$
in $\Delta$ to a family $\Phi_t(\zeta)$ of analytic functions in the unit discs, where parameter $t$ runs a real 
manifold $M^k.$
\end{Remark}
 
\section{Lower bounds for $CR$-dimensions}

In this section we apply Theorem 1 to estimating from below $CR$-dimension of manifolds in terms of families of
attached analytic discs.

\subsection{Detecting  positive $CR$-dimensions, the case $\dim \Lambda \leq n$}
We introduce two definitions, addressing to conditions 1 and 2 of Theorem 1.
 
Let $\mathcal F_{M^k}$ be a family of analytic discs attached to a real manifold $\Lambda=\Lambda^d$ and parameterized 
by a manifold $M^k.$ 

\begin{Definition}\label{D:degenerate_family} We say that the family $\mathcal F_{M^k}$ is {\bf degenerate} if there is a parameterization 
$\Phi:S^1 \times M^k \mapsto \Lambda$ satisfying the condition 1 of Theorem 1. Due to Remark \ref{R:condition_1}, this means that
either $k+1 >d$ or $k+1=d$ and $\deg \ \Phi=0.$
\end{Definition}

\begin{Definition}\label{D:homo_nontriv} We say that the  family 
$\mathcal F_{M^k}$ of analytic discs, parameterized by the mapping $\Phi,$
 attached to the manifold $\Lambda$ is {\bf homologically nontrivial}, if it possesses a parameterization
$\Phi:\Delta \times M^k \mapsto X^n$ satisfying the condition 2 of Theorem 1. Namely, this means that the image $\Phi(\{0\} \times M^k)$ of the fundamental $k$-cycle is homologically nontrivial in $\Phi(\overline \Delta \times M^k),$ relatively to $\overline \Delta \times \partial M^k.$
\end{Definition}

Remind, that a cycle $C \subset Y,$ relative to a subset $Y_0 \subset Y,$ is a chain which  boundary lies in $Y_0$.
The relative cycle $C$ is not homological to zero, relatively (to $Y_0$),if it is homological in $Y$ to no chain lying in $Y_0.$

Most transparent the condition of homologically nontriviality of the family of attached analytic discs looks when $M$ is a closed manifold. Then it  means that
the $k$-cycle $c=\Phi(\{0\} \times C)$ of the "centers"' of the discs $D_t$ bounds no $(k+1)$-chain 
in the union $\cup_{t \in M} \overline D_t$ of the closed discs.

Applying Theorem 1 to parameterization $\Phi$ of families of analytic discs we obtain its geometric version:

\begin{theorem}\label{Thm2} Let $\Lambda=\Lambda^d \subset X^n$ be a real-analytic manifold of dimension $d \leq n$, admitting a real-analytic 
$k$-dimensional  regular homologically nontrivial family $\mathcal F=\{D_t\}$ of \
attached analytic discs, parameterized by the manifold $M=M^k$.
Assume that the homology group $H_{k-2}(M,\partial M)=0.$
Suppose that the family is degenerate, i.e. 
either a) $k=d$ or  b) $k=d-1$ and $\mathcal F$ has a parameterization of topological degree $0$ on $S^1 \times M^k.$
Then the manifold $\Lambda$ has the positive $CR$-dimension at any point $b \in \Lambda,$ i.e. $\Lambda$ is nowhere totally real. 
\end{theorem}

\pf
Theorem \ref{Thm2} immediately follows from Theorem 1. Indeed, the conditions a) and b) for the family $\mathcal F$
of Theorem \ref{Thm2} translate, due to Definitions \ref{D:degenerate_family} and \ref{D:homo_nontriv}  exactly as the  
conditions of Theorem 1 for the parameterization $\Phi.$

The conclusion 
$rank_{\mathbb C} \ d\Phi < d$ on $\overline \Delta \times M^k$ of Theorem 1  implies that
for each point $u=(\zeta,t) \in \partial \Delta \times M^k$ the $CR$-dimension of $\Lambda$ at the point $b=\Phi(u)$ is
$$c(b)=dim_{\mathbb R} \ T_b\Lambda + dim_{\mathbb R} \ (i T_b\Lambda)- dim_{\mathbb R} \ (T_b\Lambda + i T_b\Lambda)=
2d- 2 rank_{\mathbb C} \ d \Phi(u) > 0.$$

%
%the vectors 
%$$\partial_r\Phi(u), \partial_{\psi}\Phi(u),\partial_{t_1}\Phi(u),\cdots,\partial_{t_k}\Phi(u)$$
%are linearly dependent over the field $\mathbb R.$ Here 
%$\zeta=re^{i\psi}$ and 
%$t_j$ 

%Taking into account the $CR$-equation $\partial_r \Phi(u)=-i\partial_{\psi}\Phi(u)$ we obtain 
%from the $\mathbb R$-linear dependence of the above system that
%the vectors 
%$$\partial_{\psi},\partial_{t_1}\Phi(u), \cdots,\partial_{t_k}\Phi(u)$$
%are linearly dependent over $\mathbb C$. But these vectors span the tangent space $T_b \Lambda$ at 
%the point $b=\Phi(u)$ and therefore
%the $d$-dimensional real space $T_b\Lambda$ is contained in a $(d-1)$-dimensional complex linear space.
%Hence
%$$2d-2c(b)=\dim_{\mathbb R} T_b \Lambda^d + 
%\dim_{\mathbb R} (iT_b \Lambda^d)-
%\dim_{\mathbb R}(T_b \Lambda \cap iT_b \Lambda)=
%2 \dim_{\mathbb C} (T_b \Lambda + iT_b\Lambda)  \leq 2(d-1)$$ 
%and 
%thus $c(b) \geq 1.$

\qed

\medskip

\noindent
\begin{Remark} \label{R:after_Thm2}

\noindent
1. One can see from the formulation of Theorem 1, that Theorem 2 which is a corollary of Theorem 1, 
remains true if we assume that $M$ and $\Lambda$ are compact 
real-analytic chains rather than manifolds. In this case  we are 
talking about $CR$-dimension at smooth points of $\Lambda.$

\noindent
2. Examples 1-4  show that the conditions of homologically nontriviality and the condition for the 
Brouwer degree in the case $k=d-1$ are essential, meaning that Theorems 1 and 2  fail to be true if any of
these conditions is omitted. 

\noindent
3. The condition  for the Brouwer degree  in Theorem 2 holds, for instance, if the
parameterizing manifold $M$ is closed, while the manifold $\Lambda$, the opposite, is not, i.e. $\partial \Lambda \neq \emptyset,$
as in Example 2. 
\end{Remark}

\medskip

In the simplest case $n=d=2$ Theorem 2 characterizes complex curves in $\mathbb C^2.$ 
Indeed, the conclusion of Theorems 1  in this case is that the tangent spaces $T_b\Lambda, \ b \in \Lambda,$  contain complex lines.
Since the real dimension of the tangent spaces is 2 this implies  that $T_b\Lambda$ is a complex line 
and therefore $\Lambda$ is a 1-dimensional complex manifold.

Thus, we have
\begin{Corollary}\label{C:Morera_curves}(Morera theorem for complex curves in $\mathbb C^2$, \cite {A1},\cite {A2}). 
Let $\Lambda $ be a real-analytic compact 2-manifold in $\mathbb C^2$ with nonempty boundary. Suppose that $\Lambda$ can be
covered by the boundaries $\gamma_t=\partial D_t$ of analytic discs $D_t \subset \mathbb C^2$, constituting a real-analytic regular
family depending on the parameter $t$ running either a curve $M^1$ or a two-dimensional manifold $M^2$.   
Suppose that in the case when 
$t$ is one-dimensional, the family $D_t$ is homologically nontrivial.
Then $\Lambda$ is a 1-dimensional complex 
submanifold of $\mathbb C^2.$
\end{Corollary}

\begin{Remark} Since the curves $\gamma_t$ in Corollary \ref{C:Morera_curves} bound the complex manifolds, $D_t$,
 for any holomorphic 1-form $\omega$ in $\mathbb C^2$ the moment condition is fullfiled
$$\int\limits_{\gamma_t}\omega=0$$
and hence Corollary \ref{C:Morera_curves}  can be viewed as Morera type characterization of complex submanifolds in $\mathbb C^2.$
\end{Remark}

The compact manifold $\Lambda$ is expected to be complex and therefore should be assumed having nonempty boundary.  
Most interesting is the case $k=1$, when the family of discs is 1-dimensional and is parameterized by a curve $M$. 
The condition b) is fullfiled, for instance,  if the family $\{D_t\}$ is closed, i.e. $M$ diffeomorphic to the circle $S^1.$ 

\medskip
\subsection{Lower bounds for $CR$-dimensions for the general case}

In the previous subsection we have considered the case $d \leq n$ and gave conditions for the $CR$-dimension of the manifold $\Lambda$ to be positive, i.e.
to be at least 1.
If $d > n$ then the $CR$-dimension is always positive and is at least $d-n$ and the lower bound $d-n$ is achieved for generic manifolds.
Therefore the lowest nontrivial bound is $d-n+1.$ 

In this subsection we consider the case $ d \geq n$ and derive  this lower estimate,
$$\dim_{CR}\Lambda \ge d-n+1,$$ 
from  the properties of  families of attached analytic discs which $\Lambda$ admits.

\medskip

Let again $$\Phi:\Delta \times M^k \mapsto \Lambda^d, \ \ k \geq d-1,$$ be a 
real-analytic regular parameterization of the family of analytic discs $\{D_t\}_{t \in M^k}$ attached to the 
real-analytic manifold $\Lambda^d \subset X^n.$

\medskip

\begin{Definition} \label{D:singular} Let $\nu \leq k$ be an integer and let $C \subset M^k$ be a chain 
in $M^k$, of the dimension $\dim \ C=\nu$ or $\dim \ C =\nu-1.$ 
We say that that the family 
$$\mathcal F_{C}=\{D_t\}_{t \in C}$$
is a {\bf singular} 
$\nu$-chain of attached analytic discs if $\Phi$ is degenerate on $S^1 \times C$ (see Definition 2.2)
and is homologically nontrivial (see Definition 2.3) on $\overline \Delta \times C.$
\end{Definition}
\medskip

The meaning of the definition \ref{D:singular} is that either $\Phi$ decreases the dimension of $S^1 \times C$ from $\nu+1$ to $\nu$,
or, when $\dim \ C=\nu-1$, there is topological 
degeneracy: $\deg \Phi=0$ on $S^1 \times C.$

\begin{Definition} \label{D:tangent} Let $\mathcal F_{C}$ be a $\nu$-chain of analytic disc from the previous definition. Denote
$\Lambda_C=\Phi(S^1 \times C)$. 
By the {\bf tangent plane} of $\mathcal F_C$
we understand 
$$Tan_b \mathcal F_{C}=T_b\Lambda_C,$$ if $b$ is a smoothness point, and  union of the tangent planes, if
$b$ is a point of self-intersection for $\Lambda_C$. 
\end{Definition}

\begin{Definition} \label{passes} Let $\Pi \subset T_b\Lambda^d$ be a tangent plane
of dimension $\nu \leq d.$
We say that the $\nu$-chain $\mathcal F_{C}$ 
{\bf passes} through the point $b$ in the  direction $\Pi$ if $\Phi \subset Tan_b \mathcal F_{C}.$
\end{Definition}

The geometric meaning of the definition is as follows. Consider all curves $\partial D_t, \ t \in C,$
passing through the point $b \in \Lambda.$  The tangent, at the point $b$, vectors to these curves are tangent
to $\Lambda.$ The same is true for the velocity  vectors (derivatives in the parameters $t_j).$  
The condition is that all these tangent vectors span the given linear subspace  $\{Pi_{\nu}$ in $T_b\Lambda.$ 

\begin{Definition} \label{admissible} We will call a real $\nu$-plane $\Pi \subset T_b\Lambda$  {\bf admissible} if there is a
complex linear subspace $P \subset T_B X^n$ such that $\Pi=P \cap T_b\Lambda$ and the intersection is transversal.
The complex dimension of such space $P$ must be $\dim_{\mathbb C}P=2n-d+\nu.$ 
\end{Definition}
\medskip

%The following example, of $U(n)$-invariant families, 
%is  pretty illustrative for  the above definitions.

\noindent
\medskip
The following theorem applies to the case $d >n$ and gives conditions to ensure that the $CR$-dimension is  bigger than the generic
dimension $d-n.$ 

\begin{theorem} Let $ d\ge n$ and let $\Lambda=\Lambda^d \subset X^n$ be a real-analytic  manifold. Suppose that $\Lambda$ 
admits a  real-analytic regular $k$-parametric, $k \ge d$, family $\mathcal F=\mathcal F_{M^k}$ of attached analytic discs.

Suppose that the family $\mathcal F$ has the following property: 
for any point  $b \in \Lambda$  and for any 
$(2n-d)$- dimensional admissible plane $\Pi \subset T_b(\Lambda)$
there exists a singular $(2n-d)$-dimensional chain
 $\mathcal F_{C} \subset \mathcal F$  passing through $b$ 
in the direction $\Pi$. Then the $CR$-dimension satisfies
$$c(b)=\dim_{\mathbb C} T^{\mathbb C}_b \Lambda \ge d-n+1$$
at any point. In other words, $\Lambda$ is nowhere generic.
\end{theorem}

\medskip

The proof of Theorem 3 follows from Theorem 2:

\pf
Suppose that $\Lambda$ is generic at  a point $b \in \Lambda.$ Then the $CR$-dimension $c(b)=d-n$ and 
the $d$-dimensional tangent plane
decomposes into the direct sum of a complex plane of the complex dimension $d-n$ and a real plane of the real dimension $d-2(d-n)=2d-n$:
$$T_b \Lambda= \Pi_{\mathbb C}^{d-n} \oplus \Pi_{\mathbb R}^{2n-d}.$$ The totally real plane $\Pi_{\mathbb R}^{2n-d}$ is obtained by complex section 
of $T_b\Lambda,$ i.e. is admissible.

By the condition there exists a singular $\nu$-chain $\mathcal F_{C}$  passing $b$ in  the direction
$\Pi^{2n-d}.$

Consider the manifold $\Lambda^{2n-d}=\Phi(S^1 \times C)$ and apply Theorem 2 to the family 
$\mathcal F_{C}$ of the analytic discs attached to $\Lambda^{2n-d}.$ 

First of all, observe that we are in the range
of dimensions required in Theorem 2, because $d \geq n$ and hence the real dimension of $\Lambda^{2n-d}$ is less than
the complex dimension of the ambient space: $2n-d \leq n$. All other conditions of Theorem 2 hold as $\mathcal F_{C}$
is a singular $\nu$-chain. 
By Theorem 2 the $CR$-dimension of $\Lambda$ satisfies $c_{\Lambda}(b) \ge 1.$ 

This means that the tangent plane $T_b\Lambda^{2n-d}$ contains
a complex line, $L$. But this tangent plane is contained in the totally real $(2n-d)$-plane 
$ \Pi_{\mathbb R}^{2n-d}$ which is 
free of nontrivial complex subspaces. This contradiction shows that $\Lambda$ is not generic at any its point. 

\qed

\begin{Remark} In the case $d=n$ Theorems 2 and 3 coincide. Indeed, then
$2n-d=d$  and therefore the only $(2n-d)$- direction is the whole tangent plane $T_b \Lambda.$ 
Also, the condition of degeneracy  in this case holds automatically with $C=M$ due to conditions a) and b).
\end{Remark}

\noindent
2. Due to real analyticity, it suffices to require that the condition in Theorem 3 holds for $b$ from an open set $U \subset \Lambda.$

\medskip

%{\bf Example 6.} Take real-analytic closed curve $$\gamma \subset \mathbb C^n$$ which is a boundary of
%an analytic disc $\gamma =\partial D, \ D \subset \mathbb C^n$ and consider the family $\mathcal F$ of analytic discs
%obtained from $D$ by linear transformation from the unitary group $U(n)$: 
%$$\mathcal F=\{u(D)\}_{u \in U(n)}.$$

\medskip
The next result generalizes Theorem 3 for higher $CR$-dimensions. 
It gives conditions for the $CR$-dimension to be at least $q$ where $q> \geq d-n$ is a given natural number not bigger than $n$.

\begin{theorem}
Let $d \geq n.$ Let the manifolds $\Lambda=\Lambda^d$ and $M^k$, and the family 
$\mathcal F_{M^k}$  be as in Theorem 3. Fix natural 
$q, \ d-n \leq q \le n.$  Suppose that for any $b \in \Lambda$
and for any admissible real plane $\Pi \subset T_b\Lambda$ of the dimension $\nu=d-2q+2,$
there exists a singular $\nu$-chain  $\mathcal F_{C} \subset \mathcal F$ 
passing through $b$ in the direction $\Pi.$ Then  the estimate for the $CR$-dimension holds:
$$ c(b) \ge q, \  b \in \Lambda.$$
\end{theorem}

\pf
Let $b \in \Lambda$. Let $c(b)$ be the $CR$-dimension at the point $b$. We know that always $c(b) \geq d-n$.
The tangent plane decomposes as in the previous theorem:
$$ T_b\Lambda=\Pi^{c(b)}_{\mathbb C} \oplus \Pi_{\mathbb R}^{d-2c(b)}$$ into a complex subspace and a totally real admissible subspace.

Applying the same argument as in the previous proof, we conclude, from the condition and Theorem 2, that  any
admissible plane $$\Pi^{\nu} \subset T_b \Lambda$$ 
of the dimension $\nu=d-2q+2,$  contains a nonzero complex subspace. 

Now suppose that, contrary to the assertion, $$c(b) \leq q-1.$$ 
Then we have $$d-2c(b) \geq d- 2(q-1)=\nu$$ 
and therefore the second, totally real, subspace, $\Pi_{\mathbb R}^{d-2c(b)},$ in the decomposition
must contain complex lines, which is not the case. This contradiction shows that $c(b) \ge q.$ 

\qed

\medskip
\subsection{Characterization of complex submanifolds and their boundaries}

Let us turn now to the case of maximally complex manifolds when the $CR$-dimension is maximally possible. If the dimension $d =\dim \Lambda$ is even, $d=2p,$ then  the $CR$-dimension $q$ takes the maximal value $q=p$ if and only if $\Lambda$ is a complex submanifold of the ambient complex manifold $X^n.$
In this case the parameter $\nu$ in Theorem 4  equals $\nu=d-2q+2=2,$ and  Theorem 4 gives the following characterization of complex submanifolds: 

 \begin{theorem} Let $\Lambda=\Lambda^{2p}$ be a real-analytic manifold, admitting a real-analytic regular family $\mathcal F$ of attached analytic discs.
 Suppose that for any $b \in \Lambda$ and for any admissible tangent 2-plane $\Pi$ there exists a singular 2-chain 
 $\mathcal F_C \subset \mathcal F$,
 passing through $b$ at the direction $\Pi$. Then $\Lambda$ is a complex $p$-dimensional manifold.
\end{theorem}

\medskip
This theorem generalizes Theorem 2 from the dimension $n=2$ to arbitrary $n.$ Notice, that the admissible 2-planes are intersections
of $\Lambda$ with complex $(n-p+1)$-planes.

Now suppose that  $\dim \Lambda=d=2p-1$ is odd. Then the maximal value for the $CR$-dimension is $q=p-1.$ The parameter $\nu$ is
$\nu=d-2q+2=3.$ The admissible 3-planes are obtained as the sections of $\Lambda$ by complex $(n-p+2)$-planes.  
Harvey and Lawson \cite{HL} proved that maximally complex real closed manifolds are exactly those which are boundaries of complex manifolds (chains).
This means that  if $\Lambda^{2p-1}$ is maximally complex then there exists $p$-dimensional complex chain $V \subset X^n$ such that
$\Lambda=\partial V.$ Therefore the special case of Theorem 4 is

\begin{theorem} Let $\Lambda^{2p-1}$ be a real-analytic closed manifold admitting a real-analytic regular family $\{D_t\}_{t \in M^k}$
of  attached analytic discs. Suppose that for any admissible tangent 3-plane $\Pi \subset T_b(\Lambda)$ there exists a singular 
3-chain $\mathcal F_{C}$ passing  through $b$ in the direction $\Pi.$
Then $\Lambda$ bounds a holomorphic p-chain, $\Lambda^{2p-1}=\partial V^p.$
\end{theorem}
\medskip

Dolbeault and Henkin \cite{DH0}, \cite{DH} proved that the 
closed maximally complex manifold $\Lambda=\Lambda^{2p-1} \subset \mathbb CP^n $ 
 bounds a holomorphic chain,
if for any complex $(n-p+1)$-plane $P,$ transversally intersecting $\Lambda,$ the 
the (1-dimensional) intersection $\Lambda \cap P$ ,  bounds a complex 1-chain.  T.-C.Dinh \cite{D1},\cite{D2},\cite{D3}
refined the result, by eliminating the condition for $\Lambda$ to be maximally complex and by using a  
narrower family of the complex $(n-p+1)$-sections.  

The Dolbeault-Henkin-Dinh results, when the manifold $\Lambda$ is contained  $\mathbb C^n$,  
correspond to  Theorem 5  for the special case of "`linear" discs. 
The required, by Theorem 5, singular chain $\mathcal F_{C},$ passing through an arbitrary point $b,$ correspond
to complex $(n-p+1)$-subspaces of the complex $(n-p+2)$-subspace defined a given admissible direction. 
The homological nontriviality can be provided by imposing the condition that
the analytic discs (complex sections) are disjoint from a fixed compact $(n-p+1)$-linearly convex set $Y$ (see \cite {D1},\cite{D2}.)

\medskip
\noindent
\begin{Remark} 

\noindent
1. Due to real analyticity, it suffices to require in Theorems 3-6 
that the singular chains condition hold for the points $b$ from an open set $U \subset \Lambda.$

\noindent
2. Most simply Theorem 6 look in the case $n=3, d=3, p=2$. Then the condition about passing 3-directions holds automatically and
Theorem is says that a sufficient condition for the
real-analytic closed 3-dimensional surface $\Lambda$ to be a boundary $\Lambda=\partial V$ of a complex 2-chain $V$
is that $\Lambda$ admits homologically nontrivial 3-dimensional regular real-analytic family of attached analytic discs. This conditions is also necessary. The corresponding family may be constructed from the analytic discs in $V,$ attached to the boundary  $\partial V=\Lambda$
and located nearby $\partial V.$ 
\end{Remark} 
  
\bigskip

\section{Morera theorems for $CR$-functions}

\medskip
In this section we apply the above results $CR$-dimensions to the  graphs of functions. As result we obtain characterization of $CR$-functions
in terms of analytic extendibility into attached analytic discs. 

\medskip
\subsection{The case of $CR-$ dimension one}

We start with the case of generic manifolds of dimension $n+1.$ In this case the $CR$-dimension equals $(n+1)-n=1.$

The following theorem follows from  Theorem 2
applied to real manifolds which are graphs.

\begin{theorem} Let $\Omega=\Omega^d \subset X^n$ be a real-analytic generic submanifold of the degree $d=n+1.$ Let 
$\{\Omega_t\}_{t \in M^k}$
be a real-analytic regular degenerate and homologically nontrivial family of attached analytic discs, 
parameterized by a compact real-analytic $k$-manifold $M^k.$
Let $f$ be a real-analytic function on $\Lambda$ such that the restriction 
$f\vert_{\partial \Omega_t}$ analytically extends in the analytic disc
$\Omega_t$, for any $ t \in M,$ or, equivalently,
$$\int\limits_{\partial \Omega_t} f \omega=0$$ for any holomorphic 1-form in $\mathbb C^.$
Then $f$ is $CR$-function on $\Omega$, i.e. $f$ satisfies the boundary $CR$-equation everywhere on $\Omega.$
\end{theorem}

\pf
Let $\Lambda$ be the graph of the function $f$, 
$$\Lambda=graph_{\Omega}f .$$
Then $\Lambda$ is a real-analytic submanifold of $X^n \times \mathbb C$ of degree $d=n+1.$ Since $\Omega$ is generic, 
we have $\dim_{CR}\Omega=d-n=1.$
On the other hand, the guaranteed lower bound for the $CR$-dimension  of the graph  $\Lambda$ is  
$$c_{\Lambda}(b) \ge d-(n+1)=0$$
and our aim is to raise the lower bound to 1.  
Denote $Q_f$ the lifting mapping
$$Q_f:\Omega \mapsto \Lambda, \ \ Q_f(z)=(z,f(z)).$$
If $F_t$ is the analytic extension of the function $f$ into the analytic disc $\Omega_t$ then the composition mapping
$$\Phi_t = F_t \circ \Psi_t$$
defines the parameterization of the family of analytic discs $$D_t=(F_t \circ \Psi(\cdot,t))(\Delta)$$
attached to the manifold $\Lambda.$ Analytic  discs $D_t$ are the graphs 
$$D_t=graph_{\Omega_t} F_t$$
of the analytic extensions $F_t$ into analytic discs $\Omega_t.$

It can be readily checked that the conditions for the manifold $\Omega$ and for the parameterization 
$\Psi$ in Theorem 7 translates as the corresponding conditions for the 
manifold $\Lambda$ and the parameterization $\Phi$ in Theorem 2 (with the dimension $n$ replaced by $n+1.$)

Then Theorem 2 implies the estimate
$$\dim_{CR} \Lambda \geq 1.$$ 
This means that at any point $u \in \Omega$  the differential $dQ_{f}(u)$ does not decrease the $CR-$ dimension and
maps the one-dimensional complex subspace in $T_u \Omega$ to a one-dimensional complex subspace of $T_{(u,f(u))} \Lambda$. In other words,
the differential $df(u)$ is a complex linear map on $T_u^{\mathbb C}\Lambda$ and therefore $f$ satisfies the tangential CR-equation at the point $u.$ 
\qed

\begin{Remark} The result similar to Theorem 7 but for special families of attached analytic discs ("thin discs"'), was obtained by Tumanov 
\cite{T1}.
\end{Remark}

\medskip
The simplest special case $n=1$  of Theorem 7 gives an answer, for real-analytic case, of a question, known as the strip-problem.  
The detailed proof is given in the article \cite {A2} and we refer the reader to them for details and the references.
 
\begin{Corollary} (the strip-problem) \label{C:4.2} Let $\Omega$ be a compact domain in the complex plane, covered by 1-parameter regular real-analytic family 
of Jordan curves $\gamma_t$, such that no point in $\Omega$ is surrounded by all curves $\gamma_t$. Then if $f$ is a real analytic function
on $\Omega$ such that all complex moments
$$\int\limits_{\gamma_t} f(z)z^m dz=0, \ \  \forall t , k=0,1,\cdots.$$
Then $f$ is holomorphic in $\Omega.$
\end{Corollary}

Let us make a few  remarks before the proof.

\noindent
\begin {Remark} 

\noindent
1.According to our considerations,  the manifold $M,$ which parameterizes
the family $\gamma_t,$  is 1-dimensional compact manifold ($k=1$), i.e. 
$M$ can be taken either the unit circle $M=S^1$ or the segment $M=[0,1].$

\noindent
2. The regularity of the family geometrically means that for each 
non-boundary point $b \Omega \setminus \partial \Omega$ and any curve $\gamma_t$ passing through
$b$ the tangent vector to $\gamma_t$ at $b$ and the 
``velocity'' vector (the derivative in $t$) are  not collinear (the "sliding points" are located on the boundary).

\noindent
3. The condition for $f$ is equivalent to analytic extendibility of $f$ inside the domain $\Omega_t$ bounded by $\gamma_t.$
By $\gamma_t$ surrounding $b$ we understand that $b$ is inside $\gamma_t$ or $b$ belongs to $ \gamma_t$, i.e. in both cases $b \in \overline \Omega_t.$

\end{Remark}

\pf {\bf of Corollary \ref{C:4.2}}
Let $\Psi$ be a regular parameterization of the family $\Omega_t.$ 
The mapping $\Psi(\cdot,t)$ can be chosen the Riemann conformal mapping of
the unit disc $\Delta$ onto $\Omega_t$. These mappings are assumed to be chosen real-analytically depending on the parameter $t$.
The 1-dimensional manifolds of parameters can be taken either  $M^1=[0,1]$ or $M^1=S^1.$ Since the compact domain $\Omega$
lies in $\mathbb R^2$ the Brouwer degree of he mapping $\Psi: S^1 \times M^1 \mapsto \Omega \subset \mathbb R^2$ equals 0.

Now, 
$CR$-functions on $\Omega$ are  holomorphic functions and 
hence to derive Corollary \ref{C:4.2} from Theorem 7, it  only remains to check that 
the family  $\Omega_t$ is homologically nontrivial. It is rigorously done in \cite{A2}. We  will give the sketch of the argument.

If the domains $\Omega_t$ constitute homologically trivial family then the  1-cycle 
$$c=\Psi (\{0\} \times M) \subset \cup_t \overline \Omega_t$$ 
is homologically trivial in $\tilde \Omega.$  The cycle $c$ is a closed curve if $M=S^1$ and if $M=[0,1]$ then $c$ is a curve with the endpoints in $$\tilde \Omega=\overline   
\Omega_0 \cup \overline \Omega_1.$$ 

The cycle $c$ can be contracted, within $\tilde \Omega$, to a point, $b$, and then by lemma about covering homotopy 
the cycle $C=\{0\} \times M$ can be correspondingly deformed to a nontrivial cycle $C^{\prime}.$ Since the cycle $C^{\prime}$ is nontrivial, it
necessarily meets each closed disc $\overline \Delta \times \{t\}$. On the other hand,  
$\Psi(C^{\prime})=\{b\}.$ Therefore $b$ belongs to each domain $\overline \Omega_t,$ 
because $\overline \Omega_t=\Psi(\overline \Delta \times \{t\}).$  We have obtained  a contradiction with
the condition and this proves -that our family is homologically nontrivial.
Thus, the family $\Omega_t$ satisfies all the conditions of Theorem 7 and Corollary \ref{C:4.2} follows.
\qed

\medskip

\subsection{Arbitrary $CR$-dimensions.}

The version of Theorem  4 for graphs is: 
\begin{theorem} Let $\Omega=\Omega^d \subset \mathbb X^n$ be a real-analytic $CR$-manifold,  $\dim_{CR} \Lambda=q.$
Suppose that $\Omega$ is covered by the boundaries of
a regular real-analytic family $\mathcal F=\{\Omega_t\}_{t \in M^k}.$ Suppose that
for any $b \in \Omega$ and for any admissible plane $\Pi=\Pi^{d-2q+2} \subset T_b\Lambda$ there exists a singular $(d-2q+2)$-chain
$\mathcal F_{C} \subset \mathcal F$ passing through $b$ in the direction $\Phi$ 
(i.e. satisfying the conditions of Theorem 4 with the parameter $q.$)  
Let $f$ be a real-analytic function on $\Omega$ satisfying the Morera condition
$$\int\limits_{\partial \Omega_t} f \omega=0, \ \forall t \in M^k,$$
for arbitrary holomorphic 1-form $\omega$,  
that is $f$ analytically extends in each analytic disc $\Omega_t.$ Then $f$ is CR-function on $\Omega.$
\end{theorem}

\pf
Denote $$\Lambda^d=graph_{\Omega} f \subset X^n \times \mathbb C.$$ Let $F_t$ be the analytic extension of $f$ into $D_t.$
Then $\Lambda$ is covered by the boundaries of the analytic discs
$$D_t=graph_{\Omega_t} F_t \subset X^n \times \mathbb C.$$
The manifolds $
\Omega$ and $\Lambda$ are linked by the diffeomorphism
$$Q_f:\Omega \mapsto \Lambda, \ Q_f(u)=(u,f(u))$$
preserving the boundaries of the attached analytic discs.

One can readily check that  the conditions of Theorem 8, for the manifold $\Omega^d$ and for the family $\{\Omega_t\}_{t \in M^k}$ ,
imply same type conditions of Theorem 4, for  $\Lambda^d$ and  for the family  $\{D\}_{t \in M^k}$,  with the same parameter $q.$ 
By Theorem 4 the $CR$-dimension of $\Lambda$ satisfies
$$c_{\Lambda}(b) \geq q.$$
But this means that the mapping $Q_f$ does not decrease $CR$-dimension and hence $Q_f$ is $CR$-mapping. Then the function $f$ is $CR$-function
as the superposition $f=\pi_2 \circ Q_f$
of two $CR$-mappings, $Q_f$ and the projection $\pi_2(u,w)=w$. 
\qed

The special case $d=2q=2n$ leads to a test for holomorphic functions, which is generalization of
Corollary \ref{C:4.2} for arbitrary dimensions. In this case $d-2q+2=2,$ the admissible 2-directions are complex lines and  we have
\begin{Corollary}($n$-dimensional strip-problem). 
Let $\Omega$ be a domain in $\mathbb C^n$ covered by the boundaries $\partial \Omega_t$ of the analytic discs from
a regular real-analytic family $\mathcal F_M.$ Suppose that $\mathcal F_M$ contains
singular 2-chains $\mathcal F_C \subset \mathcal F$ passing through each point $b \in \Lambda$ (or, at least through each point
$b$ in an open set) in any prescribed 1-dimensional complex direction.

Let $f$ be a real-analytic function in $\overline \Omega$ and assume that 
$$\int_{\partial \Omega_t} f\omega=0$$
 for every $t \in M$ and every holomorphic 1-form.
Then $f \in Hol(\Omega).$
\end{Corollary}

Another interesting 
special case of Theorem 8 is when $\Omega$ is the boundary of a complex manifold. In this case  $d=2p-1, \ q=p-1,
\ d-2q+2= 3$.
Then we obtain the following characterization of boundary values of holomorphic functions because by Bochner-Severi theorem the boundary values coincide, in smooth case, with $CR$-functions: 
\begin{Corollary}\label{V^p} Let $V \subset  \mathbb C^n$ be a $p$-dimensional complex manifold with the real-analytic boundary $\partial V=\Omega.$
Suppose that $\Omega$ admits a real-analytic regular family $\{D\}_{t \in M^k}$ of attached analytic discs and assume that
there are 3-dimensional chains passing through $\Omega$ in any admissible 3-dimensional direction. 
If $f$ is a real-analytic function on $\Omega$ and $f$ admits analytic extension into each analytic disc $D_t, \ t \in M^k$,
then $f$ extends from $\Omega$ as a holomorphic function in $V.$
\end{Corollary}

\noindent
\begin{Remark} Notice that all the analytic discs $D_t$ belong to $V.$ One of the way to provide homological nontriviality
of this family is to require that they surround a hole, meaning that the analytic discs  
$D_t, \ t \in M^k,$  fill $V \setminus V_0$ where $V_0$ is an open subset. The family in the next result is exactly of this type.
\end{Remark}

\medskip
The following example is a special case of Corollary \ref{V^p}.

{\bf Example } ({\it Globevnik-Stout conjecture} \cite{GS}). Nagel and Rudin \cite{NR} proved that if $f$ is a continuous function on the unit sphere
in $\mathbb C^n$ having analytic extension into any complex line on the fixed distance to the origin then $f$ is the boundary value of
a holomorphic function in the unit ball. The proof essentially used harmonic analysis in the unitary group and did not extend to
non-group invariant case. 

In \cite{GS} the problem of generalization of Nagel-Rudin theorem for arbitrary domains was formulated.
T.-C.Dinh \cite {D2} proved Globevnik-Stout conjecture under assumption of non real-analyticity of the line sections.
Baracco,Tumanov and Zampieri \cite{BTZ} confirmed the conjecture  with tangent Kobayashi geodesics in the place of linear sections.
(see \cite{A2} for the extended references). 

We give here the result which contains proof of the conjecture from \cite{GS} for real-analytic functions and manifolds 
(see \cite{A1}, \cite{A2}). 
The following theorem is a consequence of Corollary \ref{V^p}.
\begin{theorem} Let $D \subset \mathbb C^n, n \ge 2,$ be a bounded domain with real analytic strictly convex boundary. Let $S \subset\subset D$ be a 
real-analytic closed hypersurface. Suppose that there is an open set $U \subset D$ such that no complex line tangent to $S$ intersects
$U$ (for instance, $S$ is convex). 
Then if $f$ is a real-analytic function on $\partial D$ 
admitting analytic extension into each complex line tangent to $S$
then $f$ is the boundary value of a function holomorphic in $D.$ 
\end{theorem}
\pf
Since $\partial D$ is strictly convex, the intersections $D \cup L$ with any complex line $L \in T^{\mathbb C}S$ is an analytic disc.
The family of these discs is parameterized by the complex tangent bundle of $S$ and is real-analytic and regular. The boundaries
$L \cap \partial D$ cover the whole $\partial D$ because $S$ is a closed hypersurface. 

All we need to check is that there are 3-dimensional degenerate homologically nontrivial subfamilies passing through each point 
$b \in \partial D$ in any prescribed admissible 3-dimensional direction $\Pi$. Indeed, consider the open set  of 
complex planes $P, \ \dim_{\mathbb C} P=2,$  such that
$$P \cap U \neq \emptyset.$$  
If $P$ meets  $\partial D$ transversally then their intersection
is a real-analytic 3-manifold. Among such planes there is a plane $P$ 
passing through $b$ and such that  $$T(P \cap \partial D)=\Pi.$$  

The intersection $P \cap S$ is 3-dimensional. There is only one complex line, $L_t$ contained in the tangent plane $T_t S$. This line
$L_t$ belongs to $P.$ 
The mapping $$S \ni t \mapsto D_t=L_t \cap D$$ defines 3-dimensional subfamily of attached analytic discs, parameterized by the hypersurface $S.$ 
The parameterization $\Phi$ can be chosen so that $\Phi(0,t)=t.$ 
\footnote{Existence of a parameterization of the form $\Phi(\zeta,t)$  for the family of complex tangent lines 
requires triviality of the complex tangent bundle $T^{\mathbb C}S.$ It is so ,for instance, if $S$ is a sphere in $\mathbb C^2$.
If this is not the case then we cut the manifold $S$ into pieces with the trivial complex tangent bundles and prove the claim separately for each
corresponding bordered portion of $\partial D.$} 

This subfamily is degenerate
because the boundaries $\partial D_t$ sweep up the manifold of the same dimension 3 as the dimension of the subfamily.

It is homologically nontrivial because the 3-cycle $\Phi(\{0\} \times S)=S$ is not homological to zero in the union 
$\cap_{t \in S} \overline D_t$  of the closed analytic discs.
Indeed,  any 4-chain $V$ bounded by $S$ must contain the "hole", $U$. But $U$ belongs to the complement 
of the discs $\overline D_t$ and therefore $V$ can not be subset of their union.

Thus, all the conditions of Corollary \ref{V^p} are fullfiled and therefore $f$ extends holomorphically inside the domain $D$.
\qed

%\noindent
%\begin{Remark}. Theorem 8 can be proved in a  a stronger form. Namely, instead of domain $D \in \mathbb C^n$ one can consider
%a (strictly convex real-analytic) domain $D$ in a complex manifold $X^{p} \subset \mathbb C^n.$ 
%The closed $(2p-1)$- surface $S \subset\subset  D \subset X^p$ 
%satisfies the same conditions as in Theorem 8. The one-dimensional extension property for $f$ is assumed to hold for the family of analytic discs
%$P \cap D,$ where $P$ is taken from a smooth variety of $(n+p-1)$-complex planes in $\mathbb C^n$ tangent to $S$.  The main 3-cycles condition
%of Theorem 7 is valid. 

%Namely, the degenerate 3-cycles are formed by the sections $D \cap V$, where $V$ is taken from a variety $\mathcal P$
%of $(n-p+2)$-complex planes passing through a fixed  open set $U$ located inside the surface $S.$  Then $\dim_{\mathbb R}(V \cap \partial D)=3$ and the %$(n-p+1)$-planes $P \in \mathcal P$ contained in $V$ constitute 3-dimensional variety of analytic discs tangent to $S.$ The homological nontriviality %of such 3-cycles is provided by existence of the  "hole", $U$, not covered by the discs.

\noindent

\bigskip
\section{Proof of Theorem 1}
In this section we prove Theorem 1 from which we have derived in the previous sections all other results. 
The proof develops the main idea from   \cite{A1}, \cite{A2}. 

Let the mapping $\Phi:\Sigma=\Delta \times M^k \mapsto X^n$
be as in Theorem 1 and suppose that all the conditions of Theorem 1 are fullfiled.
We assume that $\Phi(\zeta,t)$ is not constant in $\zeta,$ i.e. $\partial_{\zeta} \Phi$ is not identically zero. 

\begin{Lemma} \label {L:L1} The following conditions are equivalent:

\noindent
(1) $rank_{\mathbb C} \ d\Phi(u) \le d-1$ for any $u \in \Sigma$  

%\noindent
%(2) the $CR$-dimension $c_{\Lambda}(b) >0$ at each point $b \in \Lambda,$ 

\noindent
(2)  $\Phi^*\eta=0$ on $\Sigma$ for any holomorphic $d$-form $\eta$ in $X^n,$

\noindent
(3) for any holomorphic $d$-form $\eta$ in $X^n$ and for any $d-1$ 
tangential smooth (real-analytic) vector fields $Z_1,\cdots,Z_{d-1}$ on the manifold $M^k$ the following identity holds:
$$J(u)=\eta(u; \partial_{\psi}\Phi(u), (Z_1\Phi)(u),\cdots,(Z_{d-1}\Phi)(u))=0.$$
Here $\ u=(\zeta,t) \in b\Sigma=S^1 \times M^k$ and $Z_j$ are viewed as tangential differential operators acting in the variable $t \in M^k$.
\end{Lemma}
\pf

\noindent
$(1) \Leftrightarrow (2).$ 

The condition (1) means that the complex dimension of the $\mathbb C$-linear span of the space $d\Phi(T_u(\Sigma))$ is less
than $d$. It occurs if and only if any differential holomorphic $d$-form vanishes on $d\Phi(T_u(\Sigma)).$ This is exactly the condition (2).

%Let $u \in \overline \Delta \times M^k$ Let $t=(t_1,\cdots,t_k)$ be local coordinates on $M^k$ near $u$.
%Let $\eta$ be a holomorphic $d$-form in $X$. Consider the function
%$$u \mapsto \eta(\partial_{\psi}\Phi(u),\partial_{t_1}\Phi(u),\cdots,\partial_{t_k}\Phi(u)), \ u=(\zeta,t) \in b\Sigma, \ \zeta=|\zeta|e^{i\psi}.$$ 
%Since $\partial_{\psi}\Phi=i\zeta\partial_{\zeta\Phi$ and $\eta$ is holomprhic form, the above defined 
%function extends  holomorphically in the disc $|\zeta|<1$. 

\noindent
$(2) \Leftrightarrow (3)$.

This equivalence  would be very easy if in (3) one would be allowed to take   $u \in \overline \Sigma$  and
vanishing of $\eta$ is checked on any system of tangent vectors. However, in (3) $u \in b\Sigma$ and only those systems
of tangent vectors are allowed which contain the vector $\partial_{\psi}\Phi$. 
Therefore, proving the equivalence (2) and (3) requires some additional efforts.

Clearly, 2) implies 3) because the vectors 
$$\partial_{\psi}\Phi(u), Z_j\Phi(u), j=1,\cdots,n$$ are tangent to $\Lambda$.
Now let us check that 3) implies 2).

First, observe that since $\eta$ is a holomorphic differential form, 
the function $J(\zeta,t)$ in (3) extends analytically from $|\zeta|=1$ to 
$|\zeta|<1$ by replacing $\partial_{\psi}\Phi(\zeta,t)=i\zeta \partial_{\zeta} \Phi(\zeta,t).$   
Therefore $J(u)=0$ for $ u \in b\Sigma$ implies $J(u)=0$ for all $u \in \overline \Sigma$.

We will prove (2) if show that 
$$\eta(u;d\Phi_u(E_0),\cdots,d\Phi_u(E_{d-1})=0$$ for any $d$ linearly  independent (over $\mathbb C$) system 
of tangent vectors $E_j \in T_u(\Sigma).$  We consider the two cases:

\noindent
1) $d=k+1.$

Then any $\mathbb C$-linearly independent system of $d=k+1$ vectors in the $k+1$-dimensional space $T_u(\Sigma)$ is equivalent to the system
$$E_0=\partial_{\zeta}\Phi(u),Z_1\Phi(u),\cdots,Z_k \Phi(u)$$
and $\eta_u=0$ on this set of vectors due to (3). This implies (2).

\noindent 
2) $d \leq k.$

In this case one has more freedom in choosing the systems of $d$ tangent vectors $E_j$ and hence
there are $\mathbb C$-linearly independent systems  $E_0,\cdots,D_{d-1}$ which do not contain the vector 
$$\partial_{\psi}\Phi(u)=i\zeta\partial_{\zeta}\Phi(u).$$ 

The condition (3) guarantees that any holomorphic $d$-form $\eta$ vanishes
on any system of the form
$$d\Phi_u(E_0),\cdots,d\Phi_u (E_{d-1}),$$ 
where $E_j$ are vectors tangent to $\Sigma.$ The important property is that the system $\{E_j\}$ includes the vector $\partial_{\psi}\Phi(u)$. 

Vanishing of any holomorphic form $\eta$  is equivalent to $\mathbb C$-linear dependence of the system \{$d\Phi(E_j)\}.$ Therefore our task now is to extent this assertion to those tangential systems which do not contain the tangent vector corresponding to differentiation in $\psi.$ .
To this end, suppose that, on the contrary, there is a holomorphic $d$-form $\eta_0$ that does not vanish on some system $\mathcal E=\{E_1,\cdots,E_d\}$
of tangent vectors to $\Sigma$ at some point $u \in \Sigma$ (and hence by continuity near $u$). 
Suppose also that this systems is not as in (3), i.e. it does not contain the vector $\partial_{\psi}\Phi.$ 

Nonvanishing of the form $\eta_0$ implies that the system $\mathcal E$ is lineary independent
over $\mathbb C$. On the other hand, we know that if we replace any of the vectors, say $E_j$, in $\mathcal E$ by the vector
$E_0=\partial_{\psi}\Phi)u),$ we obtain $\mathbb C$-linearly dependent $d$-system. This means that the vector $E_0$ belongs to
$\mathbb C$-linear span of each system $\mathcal E \setminus E_j$. Since the full system $\mathcal E$ is $\mathbb C$-linearly independent,
the intersection of these spans consists only of the vector 0. Thus, we conclude that
$$E_0=\partial_{\psi}\Phi(u)=0.$$ 

Thus, we have proven that
 $$\eta_0(u; d\Phi_u(\mathcal E) \neq 0 $$  only when $E_0=0.$
However, the zero set of function $E_0=\partial_{\psi}\Phi(\zeta,t),$ which is analytic in $\zeta$ and not identically zero,
is nowhere dense and therefore, by continuity, $\eta_0(u;d\Phi(\mathcal E))=0.$ Contradiciton. 

We have checked that, as well as in the case $d=k+1$, in the case $d \leq k$ any holomorphic form $\eta$ vanishes on the image $d\Phi_u(\mathcal E)$
of any system $\mathcal E$ of $d$ tangent vectors in $T\Sigma$. Thus, (2) holds.

%Fix $b=\Phi(u), u \in b\Sigma.$ To prove that the skew-symmetric  polilinear $d$-form
%$\eta(b)$ vanishes identically on $T_b\Lambda$ it suffices to prove that $\eta=0$ on a basis of the space $T_b\Lambda.$
%We construct this basis as follows. The first vector will be 
%$$E_0(u)=\partial_{\psi}\Phi(u) \in T_b\Lambda.$$

%The vector $E_0(u) \neq 0$ because $\Phi(\zeta,t)$ is an immersion and therefore the vector $E_0(u)=i\partial_{\zeta}\Phi(\zeta,t)$
%is never zero  as long as $\zeta \neq 0.$ 

%Further,  the mapping $\Phi :S^1 \times M^k \mapsto \Lambda^d$ is regular, i.e. has rank
%$d$. Hence, for a fixed point $u$, one can construct  $d-1$ smooth vector fields $Z_j, j=1,\cdots,d-1$ on $M^k$ so that by adding the vector 
%$E_0(u)$ we obtain the  system
%$$\mathcal E(u)=\{E_0(u),E_1(u),\cdots,E_{d-1}(u)\},$$ in which
%the vectors $E_j(u)=Z_j\Phi(u), 1 \leq j \leq d-1,$ constitute a basis in the tangent space $T_b\Lambda, b=\Phi(u).$ By condition 3)
%we have
%$$J(u)=\eta(u; \mathcal E(u))=0$$
%and 2) follows because the system $\mathcal E(u)$ is a basis in $T_b(\Lambda).$

\qed

\medskip

\subsection{The Jacobians}\label{S:Jacobians}

Suppose that  the conclusion of Theorem 1 is not true and the mapping $\Phi$ is nondegenerate.  
Remind, that in the formulation of Theorem 1 the real dimension $d$ of the manifold $\Lambda^d$ does not exceed the cmplex dimension $n$
of the ambient complex space $X^n:$ $ d \leq n.$  

By Lemma \ref{L:L1}, the equivalence 1) $\Leftrightarrow$  3), this means that
there exists a holomorphic $d$-form $\eta$ and  $d-1$ smooth tangential vector fields $Z_1,\cdots,Z_{d-1}$ on $M^d$ such that the function
$$J(u)=\eta(u;\partial_{\psi}\Phi(u),Z_1\Phi(u),\cdots,Z_{d-1}\Phi(u))$$
does not vanish identically  for $u=(\zeta,t) \in S^1 \times M^n.$ Our goal is to arrive to contradiction.

From now on, to the end of the proof, we fix the above form $\eta$ and the vector fields $Z_j, \ j=1,\cdots,n.$
Observe that the function $J$ is naturally defined in the solid manifold $\overline \Delta \times M^n$
because $\Phi(\zeta,t)$ and its derivative in the angular variable, $\partial_{\psi}=i\zeta\partial_{\zeta}$,
are defined for $|\zeta| \leq 1$, 
while the differential operators $Z_j$ act only in the variable $t \in M^n$ and do not depend on $\zeta$.

We will call the function $J$ {\it Jacobian}. Our main assumption is that $J$ is not identically zero and the goal is to obtain a contradiction
with this assumption. 

\begin{Lemma} \label{L:Jacob}
The Jacobian $J(\zeta,t) $ has the properties:

\noindent
a) $J(\zeta,t)$ is holomorphic in $\zeta \in \overline \Delta.$

\noindent
b) $J(0,t)=0, \forall t \in M^n.$

\noindent
c) the function $J^2/|J|^2$ can be represented on $(S^1 \times M^n) \setminus J^{-1}(0)$ as
\begin{equation}\label{E:Jacob}
\frac{J^2}{|J^2|}=\sigma \circ \Phi,
\end{equation} 
for some smooth function $\sigma.$

\end{Lemma}
\pf
The property a) follows from the definition (1) of the Jacobian because $\Phi(\zeta,t)$ is holomorphic in $\zeta$ and
$\eta$ is holomorphic form. The property b) is due to vanishing the vector
$$E_0(\zeta,t)=i\zeta\partial_{\zeta}\Phi(\zeta,t)$$
when $\zeta=0.$

Let us prove the property c).  
It suffices to prove that for 
 points $u_1, u_2 \in S^1 \times M^n$ such that 
$\Phi(u_1)=\Phi(u_2)$  we have
$$\frac{J^2(u_1)}{|J(u_1)|^2}= \frac{J^2(u_2)}{|J(u_2)|^2},$$
provided $J(u_1)$ and $J(u_2)$ are different from zero.

The condition $J(u_1), J(u_2) \neq 0$  implies that each of the two systems of $d$ vectors 
$$\mathcal E^{\alpha}=\{\partial_{\psi}\Phi(u_{\alpha}), Z_1\Phi(u_{\alpha}),\cdots, Z_{d-1}\Phi(u_{\alpha})\} , \ \ \alpha=1,2,$$
are linearly independent over the field $\mathbb R.$  Since all these vectors belong to the same $d$-dimensional space $T_b(\Lambda)$,
each of the two systems  constitutes a basis in this space. 

The transition $d \times d$ matrix $\mathcal A,$ from the basis $\mathcal E^1$ to the basis $\mathcal E^2,$ is real and the Jacobians
differ by the determinant of the transition matrix:
$$J(u_1)=\det \mathcal A \cdot J(u_2).$$
Since the determinant is  a real number, we have $$(\det \mathcal A)^2= |\det \mathcal A|^2$$ and the property c) follows. 
\qed

\bigskip
\subsection{The structure of the critical set $J^{-1}(0)$}

Let us start with a   brief analysis of the structure of the critical set $J^{-1}(0).$ We want to show that this set
determines $k$-current in $\overline \Delta \times M^k$ by means of the singular form $$d\ln J.$$ 
The values of this $k$-current on a $k$-form $\omega$ will be defined by
integration of the wedge-product 
$$d\ln J \wedge \omega.$$
 
We refer the reader to \cite{A2}
for more details.

Since $J$ is real-analytic, the equation $J(\zeta,t)=0, \ (\zeta,t) \in \Sigma=\overline \Delta \times M^k$
defines on $\Sigma$ an analytic set.
For each fixed $t \in M^k$ the analytic in the closed unit disc
function $J(\zeta,t)$ either has finite number of isolated zeros in $\overline \Delta$ or vanishes  identically.

Correspondingly,  the zero set of $J$ can be decomposed into two parts:
$$J^{-1}(0)= N_1 \cup N_2,$$
where $N_1$ is the $k$-chain defined by the isolated zeros $\zeta \in \Delta$, and $N_2$ consists of the discs
where $J$ vanishes identically with respect to $\zeta$. The set $N_2$ has the form of the direct product:
$$N_2=\overline \Delta \times T$$
of the closed unit disc and an analytic set T$\subset M^k.$ 

The pieces $\overline \Delta \times T^s,$ where $T^s$ are strata (see \cite{GM}) 
of dimensions $s \leq k-3,$ are negligible for integration of $k$-forms,  because
they have the dimension $2+s < k.$  So only the strata of the dimensions $s=k-2$ and $s=k-1$ can contribute to the integration.

The case $s=k-1$ will be removable  due to the following argument. The strata $T^{k-1}$ is defined, at least locally, by zeros of a real function,
$T^{k-1}=\{\rho=0\}$ so that  $J(\zeta,t)=\rho(t) \ I(\zeta,t)$ where $I(\zeta,t)$ does not vanish identically in $\zeta$ for $t \in T^{k-1}.$
Then the the real factor $\rho$ cancels in the ratio 
$$\frac{J}{\overline J}=\frac {I}{\overline I}.$$
In our construction n the sequel we  will consider the 
normalized vector field $J^2/|J^2|=J/\overline J$  rather than $J$ so that the above type of singularity is removable.

From the point of view  of currents, removability of this type of singularity can be explained as follows. We have
$$d \ln (J/\overline J)=d\ln J^2 -d\ln |J|^2$$
and $d\ln |J|^2$ defines the zero current because the singularity of $\ln|J|^2$ is integrable in dimensions greater than 1 and hence removable. Thus,
the $(k-1)$-dimensional strata in $N_2$ are removable for the current $2d\ln J.$  
 
Thus, only the $k$-dimensional part
$$N_2=\overline \Delta \times T^{k-2}$$
may contribute in the current $d\ln J.$

\subsection {The dual cocycle}

\medskip
We will use Poincare duality between the relative, with respect to $A$, homology group $H_k(Y,A; \mathbb Z)$ and the 
homology group $H_{n-k}^{open}(Y \setminus A;\mathbb Z)$ of the complement of $A$, of the complementary dimension. The duality is provided
by the pairing: the intersection index of cycles and cocycles. 

The Poincare duality implies that a relative, to $A$, cycle $C$ is homologically nontrivial if and only if there exists
a dual cycle (cocycle) $C^{\perp}$, of the complementary dimension, disjoint from $A$ and having with $C$ the nonzero intersection index.

\begin{Lemma}\label{L:dual_cocycle} 
let $S$ be a 2-chain in the $(k+2)$-manifold $\overline \Delta \times M^k$, having the following properties:

\noindent
1) $S$ transversally intersects the ($k$-dimensional) critical set $J^{-1}(0)$, at finite number of point from $\overline \Delta \times M^k$,

\noindent
2) $\Gamma=\partial S \subset S^1 \times (M \setminus \partial M),$

\noindent
3) $\Gamma$ is the union of $\Phi$-fibers, i.e. $\Gamma=\Phi^{-1}(\Phi(\Gamma)).$

Then the intersection index $ind (J^{-1}(0) \cap S)=0.$
\end{Lemma}

\pf
First we have to explain how should one understand the intersection index at the boundary points $b \in \partial S$.
The function $J(\zeta,t)$ is holomorphic in the disc $|\zeta| <1 +\varepsilon$ and each zero $a \in \partial S$ of $J$ is a 
non-boundary isolated zero
for the holomorphic extension of  $J$. We can also extend the surface $S$ in the domain $|\zeta|<1+\varepsilon$ as a chain $\tilde S$.
Then we set
$$ind_a (J^{-1}(0) \cap S)=\frac{1}{2} \  ind_a (J^{-1} \cap \tilde S).$$
The non-boundary common points of $J^{-1}(0)$ and $S$ are exactly (isolated ) zeros of  $J$ on the surface $S.$

Now, the boundary $\Gamma$ of the surface $S$ is the trace of $S$ on the boundary of the solid manifold $\overline \Delta \times M^k$:   
$$\Gamma=[S \cap (\overline \Delta \times \partial M^k)] \cup [S \cap (S^1 \times (M^k \setminus \partial M^k)].$$
By the condition  2) the first part is empty and hence

$$\Gamma=S \cap (S^1 \times (M^k \setminus \partial M^k)).$$
This yields that is the boundary of $S$ is the finite union of closed curves:
$$\Gamma=\Gamma_1 \cup \cdots \cup \Gamma_q. $$

Hence, the intersection index $ind (J^{-1}(0) \cup S )$ 
equals the algebraic sum of the indices  $ind_ {S, a_j} (J)$ of the vector field $J$ on the surface $S$ at the isolated zeros $a_j$:
\begin{equation}\label{E:2}
ind \ (J^{-1}(0) \cap S)= \sum\limits_j ind_  {S, a_j}  (J).
\end{equation}
By definition of the index of vector field at an isolated zero, we have for the vector field $J^2$: 
\begin{equation}\label{E:3}
2 \ ind_{S, a_j} \ (J)= ind_{S, a_j} \ (J^2)=ind_{S,a_j} \ (\frac{J^2}{|J^2|}).
\end{equation}

The sum in the right hand side coincides with the winding number $W_{\Gamma}(J^2/|J^2|)$ of the normalized vector field
$J^2/|J^2$ along the closed curve $\Gamma=\partial S.$ Therefore formulas (\ref{E:2}) and (\ref{E:3}) yield
\begin{equation}\label{E:4}
2 \ ind \ (J^{-1}(0) \cap  S)=W_{\Gamma}(\frac {J^2}{|J^2|})=\frac {1}{2 \pi i} \int\limits_{\Gamma}\frac{ d(J^2/|J^2)}{J^2/|J^2|}.
\end{equation}

Recall that by Lemma \ref{Jacob} c), formula \ref{E:Jacob}, 
the function $J^2/|J^2|$ is represented on $S^1 \times M^k$, out if the critical set $J^{-1}(0)$, as the superposition
$J^2/|J^2|=\sigma \circ \Phi$ 
and, in particular,  this representation takes place on the subset $\Gamma \subset S^1 \times M^k.$

\medskip
According to Remark \ref{R:condition_1} the degeneracy condition 1 in Theorem 1 means that one of the two following cases take place. 

\noindent 
{\it The case a), \ $k=d= \dim \Lambda.$} 
The regular mapping $\Phi$ degenerates, meaning that it maps the $(d+1)$-dimensional manifold
$S^1 \times M^d$ 
to the $d$-dimensional  manifold $\Lambda^d$ and therefore the $\Phi$-fibers 
$\Phi^{-1}(b) \cap (S^1 \times M^k),$  corresponding to non-boundary  
points $b \in \Lambda^d,$  are 1-dimensional. 

But the boundary $\Gamma=\partial S \subset S^1 \times M^k$ is also 
1-dimensional and due to the condition  2)  we conclude that  $\partial S$ is union of $\Phi$-fibers belonging to $S^1 \times M^k.$ 

However, from Lemma \ref{L:Jacob}, formula (\ref{E:Jacob}), 
we know that the function $J^2/|J^2|$ is constant on these fibers and therefore this function is piece-wise constant
on $S$.
Since $J^2/|J^2|$ is continuous, it is constant on each curve $S_j \subset \partial S$ and hence has zero change of the argument on $S$.
Thus,
$$2 W_{\Gamma}(J) = W_{\Gamma}(\frac {J^2}{|J^2|})=0.$$
Hence the right hand side in (\ref{E:4}) equals $0$ and  this proves Lemma \ref{L:dual_cocycle} for the case a).

\medskip
\noindent {\it The case b), \ $k=d-1$}. Then $\Phi$ maps the $d$- dimensional manifold $S^1 \times M^{d-1}$ to the  manifold $\Lambda^d$ 
of the same dimension.
There is no degeneracy from the point of view of dimensions, but  by the condition $\Phi$ degenerates in the sense that the Brouwer degree 
of $\Phi$ on $S^1 \times M^{d-1}$ is 0, 
$$\deg \Phi \vert_{S^1 \times M^{d-1}}=0.$$
Then, again, from formula (\ref{E:4}) we have

$$2 \ ind \ (J^{-1}(0)\cap  S)= W_{\Gamma}(\frac{J^2}{|J^2|})=W_{\Gamma}(\sigma \circ \Phi)=
\deg \Phi \cdot W_{\Phi(\Gamma)}(\sigma)=0.$$
Lemma is proved.
\qed
\bigskip
\subsection { End of the proof of Theorem 1.}

\medskip
First of all, by real-analyticity, $rank \ d\Phi(u)$ either equals $k+1$ everywhere in $\overline \Delta \times M^k$ or equals $k+2$
everywhere except proper analytic subset $Sing \ \Phi.$ Therefore the image $\Phi(\overline \Delta \times M^k)$ is 
either $(k+1)$-dimensional manifold,
or it is $(k+2)$-dimensional with singularities on the set of critical values $\Phi (Sing \  \Phi).$

By the condition of homological nontriviality, the mapping $\Phi$ maps the fundamental $k$-cycle $C=\{0\} \times M^k$
onto cycle $\Phi(C)$ which represents a nonzero element in the $k$-th homology group of $G=\Phi(\overline \Delta \times M^k)$ relatively to 
$G_0=\Phi(\overline \Delta \times \partial M^k).$  

By property b) of Lemma \ref{L:Jacob}, $C$ is contained in the critical set $J^{-1}(0).$ This set consists of two parts
$$ J^{-1}(0)=N_1 \cup N_2,$$
where the first, regular, part $N_1$ contains $C$. Since $\Phi(\zeta,t)$ is holomorphic in $\zeta$
the chain $\Phi(N_1)$ is cooriented with $\Phi(C)$ and hence is not homological to 0. 

As for the second, the "zero-disc", part  $N_2=\overline \Delta \times T^{k-2}$ is concerned, it
represents the class $0$ in the relative homology group 
$$H_{k}(\overline \Delta \times M^k, \overline \Delta \times \partial M^k)$$ 
because $\overline \Delta $ is contractible and 
$$[T^{k-2}] \in H_{k-2} (M^k,\partial M^k)=0.$$ 
Then the image chain
$\Phi(N_2)$ is homological to $0$ in the corresponding homology groups of the image. 
Thus, we conclude that the entire chain $\Phi(J^{-1}(0)$ represents a nonzero $k$-homology class
$$[\Phi(J^{-1}(0)]=[\Phi(N_1)] \neq 0.$$

Now we use  the Poincare duality.  There exists a dual chain $$(\Phi(N_1))^{\perp}= R \subset \overline \Delta \times M^k$$  having the nonzero intersection index with $\Phi(N_1)$ and disjoint from $\Phi(S^1 \times \partial M).$ 
Remind that the "zero-disc" part $N_2$ contributes nothing in the intersection index as it is homologically trivial.

Thus, we have  a chain $R$ of the complementary dimension
$$\dim R=\dim \ G  - dim \ \Phi(J^{-1}(0))=\dim \ G - k$$
having the properties:

\noindent
1) $R \cap (\overline \Phi(\Delta \times \partial M^k)) =\emptyset,$

\noindent
2) $R$ transversally intersects $\Phi(N_1)$  at a a finite set of points,
 
\noindent 
3) the intersection index
\begin{equation}\label{E:5}
ind \ (R \cap \Phi(J^{-1}(0)))=ind \ (R \cap \Phi(N_1)) \geq ind \ (R \cap \Phi (\{0\} \times M^k)) > 0. 
\end{equation}

Define $S$ as the preimage of $R$:
$$S=\Phi^{-1}(R).$$
Then $S \subset \overline \Delta \times M^k.$ 
We want to understand the dimension of $S.$ 
Denote $$m=\dim \ \Phi(\Delta \times M^k).$$ 
Then the dimension of $\Phi$-fibers is $k+2-m.$
Since $R$ has the complementary dimension,
$$\dim \ R=\dim \ \Phi(\Delta \times M^k)-k=m-k,$$
its preimage has dimension 2:
$$\dim \ S=\dim \ R + (k+2-m)=(m-k)+(k+2-m)=2.$$
Thus, $S$ is a 2-dimensional chain.
By the construction, $S$ and $\Gamma=\partial S$  are unions of $\Phi$-fibers.  The conditions 1),2),3)  imply 
$$S \cap (\overline \Delta \times M^k)=\emptyset.$$
Now, $\Phi$ is regular on $S^1 \times M^k$ and therefore $\Gamma$ meets $J^{-1} \cap (S^1 \times M^k)$ transversally.
By small deformation of $S$ without changing its boundary $\Gamma$ one can make  $S$ intersecting $N_1$ transversal
at each point. Such perturbation preserves the property that $\Gamma$ is the union of $\Phi$-fibers,  i.e. $\Gamma=\Phi^{-1}(\Phi(\Gamma))$

Thus, $S$ satisfies all the conditions of Lemma \ref{L:dual_cocycle}  and hence
\begin{equation}\label{E:6}
ind \ (S \cap N_1)=0.
\end{equation}

We  have obtain from (\ref{E:5}) and  (\ref{E:6}) 
the two contradictory conclusions for  the intersection indices: $ind (\Phi(S) \cap \Phi(N_1))=ind (R \cap \Phi(N_1)) > 0$,
on one hand, 
and $ind(S \cap N_1)=0,$ on the other hand.  

This contradiction shows that our main assumption that $J$ does not vanish identically can not be true. Therefore $J(u) =0$
for all $u \in \Delta \times M.$ Due to Lemma \ref{L:L1} 
and the construction of the Jacobian $J$ in Lemma \ref{L:L1}, 3), this completes the proof of Theorem 1.
\begin{figure}[h]
	\centering
\scalebox{0.5}{		\includegraphics{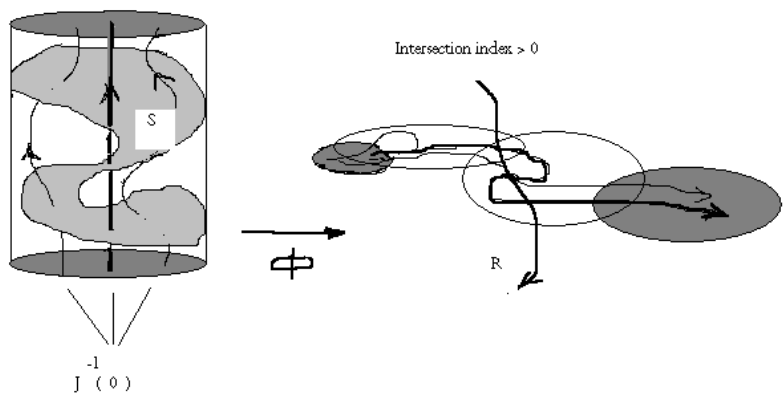}}
	\caption{The critical set $J^{-1}(0)$ and the chains $R$ and $S$.}
	\label{fig:fig.2}
\end{figure}

\section{Appendix: An alternative proof of Theorem 1 for the case $\partial M=\emptyset$}
For the case of closed families, "more analytic"' proof can be given, using linking numbers and de Rham duality  
rather than  intersection indices and Poincare duality.

\begin{Lemma}\label{L:stokes}
Suppose that $\partial M=\emptyset.$ 
For any closed differential $k$-form $\omega$ in a neighborhood of  $\Phi(\overline \Delta \times M^k).$
Then
\begin{equation}\label{E:7}
\int\limits_{J^{-1}(0)} \Phi^*\omega=0.
\end{equation}
In particular, the chain $J^{-1}(0)$ is a cycle.
\end{Lemma}
\pf
The last assertion follows form de Rham duality, as a $k$-chain defines the integer valued current on 
closed differential $k$-forms if and only if it is a cycle.

Now, the Stokes formula for the currents yields
\begin{equation} \label{E:8}
0=\int\limits_{(\overline \Delta \times M^k)\setminus J^{-1}(0)} d(d\ln J^2 \wedge \Phi^*\omega)=
\int\limits_{S^1 \times M^k} d\ln J^2 \wedge \Phi^*\omega + 2\int\limits_{J^{-1}(0)} \Phi^*\omega.
\end{equation}
It remains to prove that the first integral in the last expression is 0.
Since the forms 
$$d\ln J^2=2\frac{dJ}{J}$$
and 
$$d \ln \frac{J^2}{|J^2|}=\frac{dJ}{J} - \frac {d\overline J}{\overline J}$$
define the same currents ,
the integral that we are interested in equals to
\begin{equation}\label{E:9}
 I=\int\limits_{S^1 \times M^k} d\ln\frac{J^2}{|J^2|} \wedge \Phi^*\omega.
\end{equation}
Using the representation (\ref{E:Jacob}) 
for $J^2/|J^2|$ from Lemma \ref{L:Jacob},  we obtain that the differential form under the sign of the integral
represents as 
\begin{equation}\label{E:10}
\Xi=\Phi^*( d\ln\sigma \wedge \omega)). 
\end{equation}

Suppose that we are in the situation of Theorem 1, that is $k=d.$
Then the degree  of the form $\omega$ is $k=d$ and therefore
$$\deg (d\ln\sigma \wedge \omega)=d+1.$$ Since the dimension of the manifold
$\dim \Phi(S^1 \times M^d)=\dim \Lambda^d=d$ is less than the degree $d+1$ of the form, the restriction of the form is zero:
$$d\ln\sigma \wedge \omega\vert_{\Phi(S^1 \times M^d)}=0$$ 
and then the form $\Xi=0$ and the integral $I=0.$

\medskip
Now suppose that the case b) take place, i.e.
$k=d-1$ and $\deg \Phi=0$ on $S^1 \times M^k.$

In this case the integral $I$ in (\ref{E:9}) is zero by the different reason, namely, because the Brouwer degree of the mapping 
$$\Phi: S^1 \times M^{d-1} \mapsto \Lambda ^d, $$
of two $d$-dimensional mappings is zero.

Indeed, from formula (\ref{E:4}) we have
$$I=\int\limits_{(S^1 \times M^{d-1}) \setminus J^{-1}(0)} \Xi= 
\deg \Phi \ \int\limits_{\Phi((S^1 \times M^{d-1}) \setminus J^{-1}(0))} d\ln \sigma \wedge \omega=0,$$ 
because $\deg \Phi=0.$ Lemma is proved.
\qed

\medskip
The final part of the proof is as follows.
By the condition the image  $\Phi(\{0\} \times M^k)$ of the fundamental cycle is homologically nontrivial
in $G=\Phi(\overline \Delta \times M^k).$ 

By de Rham duality there exists a closed nonexact $k$-form $\omega$in a neighborhood of $g$
such that
$$\int\limits_{\Phi(\{0\} \times M^k)} \omega > 0.$$
As above, the cycle of integration consists of 2 parts:
$$\Phi(J^{-1}(0))=\Phi(N_1) \cup \Phi(N_2).$$ 
The first part is a cycle cooriented with $\Phi(\{0\} \times M^k)$ because  $\Phi(\zeta,t)$ is holomorphic in $\zeta.$
The second part $\Phi(N_2)$ is negligible because is represents the zero homology class.
Therefore, we have
\begin{equation}\label{E:11}
\int\limits_{J^{-1}(0)}\Phi^*\omega=
\int\limits_{\Phi(J^{-1}(0))}\omega  >0.
\end{equation}

This contradiction with Lemma \ref{L:stokes} shows that the assumption that $J$ is not identically zero is not true and this completes the proof,
due to Lemma \ref{L:Jacob} and the definition of the Jacobian $J$ in Section \ref{S:Jacobians}

\medskip
\section{Concluding remarks.}

In Theorem 1, we consider the $(k+2)$-dimensional manifold $\Sigma=\Delta \times M^k$ which carries a natural $CR$- structure. 
Essentially, what we prove in Theorem 1 is that if a $CR$-mapping $\Phi$ is not degenerate (has the maximal rank) on
$\Sigma$   and if $\Phi$ induces nontrivial   homomorphism of the  homology groups $H_k(\overline \Sigma) \cong H_k(M^k)$ to
the group $H_k(\Phi(\overline \Sigma))$ 
(i.e. $\Phi$ is homologically nontrivial), then for the image of the boundary $b\Sigma=S^1 \times M$ 
we have $H_{k+1}(\Phi(b \Sigma)) \neq 0.$ 

Theorem 1 deals with  $CR$-manifolds of the simple form, namely with those  manifolds which are 
the  direct products  $\Sigma=\Delta \times M$ 
of the unit disc and a real manifold of parameters. We think that it can be generalized to more general $CR$-manifolds. 
The expected theorem would state that, under certain conditions,
nondegenerate $CR$-mapping $\Phi$ of a bordered $CR$-manifold can not be homologically trivial on the boundary. 

According to our discussion in Introduction, the expected result, translated to  the language of the $\Phi$-images, 
would lead to a generalization of the weak form of the argument principle. This generalization
would claim that the boundary of a nontrivial, in a certain sense, $CR$-manifold can not be homologically trivial. 

Moreover,  we conjecture that a strong, quantitative, version of the argument principle for $CR$-mappings and $CR$-manifolds
might be realized.  If in the case of  complex manifolds,
the argument principle can be interpreted as the equality between the linking number of the
boundary of a complex manifold with respect to another 
complex manifolds, on one hand, and the 
number of the intersection points of the two complex manifolds, on the other hand, 
then in  more general case of $CR$-manifolds, one may expect a combined relation.
In this conjectured relation, the complex and real counterparts of the manifold, coming correspondingly from  complex and totally real tangent bundles,  
would contribute  in a different way.  Such generalization is supposed to be the subject of a future publication.

\medskip
{\it Ackhowledgements.} I am grateful to V. Gichev, J. Globevnik, T.-C. Dinh, P. Dolbeault, L. Ehrenpreis, G. Henkin,
P. Kuchment, C. Laurent, S. Shnider, E. ~L. Stout, A. Tumanov, M. Zaidenberg and L. Zalcman for their interest 
to this work and useful discussions and remarks. 
%
% for drawing my attention to his articles \cite{D1},\cite{D2},\cite{D3} and
%Prof. G. Henkin and Prof. P. Dolbeault for the fruitful discussion of
%the results of their articles \cite{DH0},\cite{DH} in connection with the subject of this article. I thank Professors
%A. Boggess, P. Kuchment, C. Laurent, S. Shnider and M. Zaidenberg for their interest to the work.
%I am grateful to Josip Globevnik for an useful remark.

\noindent
e-mail:agranovs@math.biu.ac.il


\begin{thebibliography}{15}

\bibitem{A1} M. Agranovsky, {\em $CR$ foliations, the strip-problem
and Globevnik-Stout conjecture}, C. R. Acad. Sci. Paris, Ser.I 343 (2006),
91-94.
\bibitem{A2} M.Agranovsky, {\em Propagation of boundary $CR$-foliations and Morera type theorems for manifolds
with attached analytic discs}, Advances in Math., 211, 1, (2007), \ 284-326.
\bibitem{AV} M. L. Agranovsky and R. E. Val'sky, {\em Maximality of invariant
algebras of functions}, Siberian Math.J.,  12 (1971),  1-7.
\bibitem{ABC} M. L. Agranovsky, C. Berenstein and D.-C. Chang, 
{\em Morera theorem for holomorphic $H^p$ spaces in the Heisenberg group},
J.Reine Angew.Math.,443(1993), 49-89.
\bibitem{AN} M. L. Agranovsky and E.K. Narayanan, {\em Isotopic families of
contact manifolds for elliptic PDE}, Proc. AMS,134( 2006), 2117-2123.
\bibitem{AS} M. L. Agranovsky and A. M. Semenov, {\em Holomorphy on unitary-
invariant families of curves in $\mathbb C^n$}, Siberian Math.J., \ 29 (1988),
149-152.
\bibitem{MLA} M. L. Agranovsky, {\em Remarks about one-dimensional holomorphic extension},
talk at Conference ``Complex Analysis and Mathematical Physics'', Divnogorsk,
1987; a note (jointly with A.M. Semenov) in Proceedings of the same conference, published by Institute
of Physics, Krasnoyark, 1987,16-18.
\bibitem{AS1} M. L. Agranovsky and A.M. Semenov, {\em Boundary analogues
of the Hartogs' theorem}, Siberian Math.J.,  32 (1991), no.1 \ 137-139.
\bibitem{AG} M.L. Agranovsky and J.Globevnik, {\em Analyticity on circles
for rational and real-analytic functions of two real variables},
J. Analyse Math. 91 (2003),\ 31-65.
%\bibitem{AYu} L. A. Aizenberg and A. P. Yuzhakov, {\em Integral representations
%and residues in multidimensional comlex analysis}, Amer.Math.Soc.,
%Providence, R. I., 1983.
%\bibitem{AH} R. A. Ayrapetyan and G. M. Henkin, {\em Analytic
%continuation of CR functions through the ``edge of the wedge''},
%Dokl. Akad. Nauk SSSR 259 (1981), 777-781.
\bibitem{AW} H.Alexander and J.Wermer, {\em Linking numbers and
 boundaries of varieties}, Ann. of Math (2) 151
(2000), 125-150.
\bibitem{BTZ} L.Baracco, A. Tumanov, and G. Zampieri, {\em
Extremal discs and the holomorphic extension from convex hypersurfaces}, Ark. Mat. 45 (2007),1-13.
%\bibitem{BW} A. Browder, J. Wermer, {\em Some algebras of functions on an arc}, J. Math. Mech. 12 (1963), 119-130.
\bibitem{D1} T.-C. Dinh, {\em Probl\'eme du bord dans l'espace projectif complexe}, Ann. Inst. Fourier, Grenoble, 48,5 (1998), 1483-1512. 
\bibitem{D2} T.-C.Dinh, {\em Conjecture de Globevnik-Stout et theoreme
de Morera pur une chaine holomorphe}. Ann. Fac. Sci. Toulouse Math. (6)
8 (1999), no.2, 235-257.
\bibitem{D3} T.-C. Dinh, {\em Sur la caracterisation du bord d'une chaine holomorphe dans l'espace projectif},
Bull. Soc. Math. de France, 127 (1999), 519-539.
\bibitem{DH0}P. Dolbeault and G. Henkin, {\em Surfaces de Riemann de bord donn\'e dans $\mathbb C \mathbb P^n$},
Aspects of Mathematics, 26, Viewweg (1994), 163-187.
\bibitem {DH} P. Dolbeault and G. Henkin, {\em Chaines holomorphes de bord donn \'e dans $\mathbb C \mathbb P^n,$}
Bull. Soc. Math. de France, 125 (1997), 383-445.
\bibitem{E2} L.Ehrenpreis, {\em Three problems at Mount Holyoke},
Contemp. Math. 278 (2001), 123-130.
\bibitem{E1} L. Ehrenpreis, {\em The Universality of the Radon Transform},
Oxford Univ. Press,2003.
\bibitem{F} F. Federer, {\em Geometric Measure Theory}, Grundlehren der Math. Wiss., 285, Springer, Berlin-
Heidelberg-New York, 1988.
\bibitem{G1} J. Globevnik,{\em Analyticity on rotation invariant families
of circles},
Trans. Amer. Math. Soc. 280 (1983), 247-254.
\bibitem {G2} J. Globevnik,{\em A family of lines for testing holomorphy
in the ball of $\mathbb C^2$}, Indiana Univ. Math. J. 36 (1987), no.3,639-644.
Trans. Amer. Math. Soc. 280 (1983), 247-254.
\bibitem{JG} J. Globevnik, Talk at Analysis Seminar, Bar-Ilan University,1987.
\bibitem{G3} J.Globevnik, {\em Testing analyticity on rotation invariant
families of curves}, Trans. Amer. Math. Soc. 306 (1988), 401-410.
\bibitem{G4} J. Globevnik, {\em Holomorphic extensions and
rotation invariance}, Complex Variables, 24 (1993), 49-51.
\bibitem{G5} J. Globevnik, {\em A boundary Morera theorem},
J. Geom.Anal. 3(1993), no.3,269-277.
\bibitem{G6} J. Globevnik, {\em Holomorphic extensions from open families
of circles}, Trans. Amer. Math. Soc. 355 (2003),  1921-1931.
\bibitem{G7} J.Globevnik, {\em Analyticity on translates of Jordan  curves},
preprint, ArXiv, math.CV/0506282; Trans. Amer. Math. Soc. (to appear).
\bibitem{GS} J.Globevnik and E.L.Stout, {\em Boundary Morera theorems for
holomorphic functions of several complex variables},
Duke Math. J. 64(1991), no 3, 571-615.
\bibitem{GS1} J. Globevnik and E. L. Stout, {\em Discs and the Morera 
property}, Pacific J. Math., 192(2000), 65-91.
\bibitem{GM} M.Goresky, R. MacPherson, {\em Stratified Morse Theory},
Springer Verlag, Berlin, Heidelberg, 1988.
\bibitem{Ha} R. Harvey,{\em Holomorphic chains and their boundaries}, Proc. Symp. Pure Math., 30, vol. 1 (1977), 309-382.
\bibitem{HL} R.Harvey and B.Lawson, {\em On boundaries of complex analytic varieties I}, Ann. of Math., 102(1975), 233-290. 
\bibitem{KM} A.M. Kytmanov and S.G. Myslivets, {\em Higher-dimensional boundary analogs of the Morera theorem in problems of analytic continuation
of functions}, Complex analysis, J.Math.Sci (N.Y.),
120 (2004), 1842-1867.
%\bibitem{LC}C.H. Look and T.D. Zhong, {\em An extension of Privalov
%theorem}, Acta Math. Sinica, 7 (1957),144-165.
%\bibitem{MS} J.W. Milnor, J.B.Stasheff, {\em Characteristic classes},
%Annals of Math. Studies, Princeton Univ. Press and Univ. of Tokyo 
%Press, Princeton,NJ,1974.
\bibitem{NR} A. Nagel and W. Rudin, {\em Moebius-invariant function
spaces on balls and spheres}, Duke Math. J. 43 (1976), 841-865.
\bibitem{R} W. Rudin, {\em Function theory in the unit ball of $\mathbb C^n$},
Springer-Verlag, Berlin, Heidelberg, New York, 1980.
\bibitem{St2} E.L. Stout,{\em The boundary values of holomorphic functions
of several complex variables}, Duke Math.J. 44 (1977), no.1, 105-108.
\bibitem{St1} E.L. Stout,{\em Boundary Values and Mapping Degree},
Michigan Math. J. 47 (2000), 353-368.
\bibitem{Sp} E. Spanier,{\em Algebraic topology}, McGraw-Hill, New York,
1966.
%\bibitem{Sh} H. Shapiro, {The Schwarz function and its generalizations
%to higher dimensions}, Wiley,1992.
\bibitem{T1} A.Tumanov,{\em A Morera type theorem in the strip},
Math.Res. Lett., 11 (2004), no. 1, 23-2
\bibitem{T2} A.Tumanov, {\em Testing analyticity on circles},  Amer. J. of Math. 129 (2007),785-790. 
\bibitem{T3} A.Tumanov, {\em Thin discs and a Morera theorem
for CR-functions}, Math.Z.,226 (1997), no. 2, 327-334.
%\bibitem{Hu} H.Sze-Tsen, {\em Homotopy theory},Academic Press, New York and
%London,1959.
\bibitem{W} J. Wermer, {\em The hull of a curve in $\mathbb C^n$}, Ann. of Math., 68 (1958), 550-561. 
\bibitem{Z1} L. Zalcman,{\em Analyticity and the Pompeiu problem},
Arch. Rat. Mech. Anal. 47 (1972), 237-254.

\end{thebibliography}
\end{document}